\documentclass{article}
\usepackage{sw20bams}
\usepackage{latexsym}
\usepackage{amsmath}
\usepackage{graphicx}
\usepackage{amsfonts}
\usepackage{amssymb}
\newtheorem{theorem}{Theorem}
\newtheorem{acknowledgement}[theorem]{Acknowledgement}

\newtheorem{corollary}[theorem]{Corollary}

\newtheorem{definition}[theorem]{Definition}

\newtheorem{lemma}[theorem]{Lemma}
\newtheorem{notation}[theorem]{Notation}

\newtheorem{remark}[theorem]{Remark}

\begin{document}

\author{Andrey Todorov\\University of California\\Department of Mathematics,\\Santa Cruz, CA 995064\\Bulgarian Academy of Sciences\\Institute of Mathematics\\Sofia, Bulgaria}
\title{Ray Singer Analytic Torsion of Calabi Yau manifolds I.}
\maketitle
\begin{abstract}
In this paper we generalized the variational formulas for the determinants of
the Laplacians on functions of CY metrics to forms of type (0,q) on CY
manifolds. We also computed the Ray Singer Analytic torsion on CY manifolds we
proved that it is bounded by a constant. In case of even dimensional CY
manifolds the Ray Singer Analytic torsion is zero. The interesting case is the
odd dimensional one.
\end{abstract}
\tableofcontents

\section{Introduction.}

One of the most remarkable formula that I encounter is the Kronecker limit
formula. It states that if

\begin{center}
$E(s)=\underset{n,m\in\mathbb{Z}}{\sum^{\prime}}\frac{1}{|n+m\tau|^{2s}}$
\end{center}

where $\tau\in\mathbb{C}$, $\operatorname{Im}\tau>0$ and ' means that the sum
is taken over all pair of integers $(m,n)\neq(0,0),$ then $E(s)$ has a
meromorphic continuation in $\mathbb{C}$ with only one pole at $s=1$ and
$\exp(-\frac{d}{ds}E(s)|_{s=0})=\left(  \operatorname{Im}\tau\right)
^{2}\left|  \eta\right|  ^{4}$ where $\eta$ is the Dedekind eta function.

It is a well know fact that in the case of elliptic curves $\{E_{\tau
}=\mathbb{C}/(n+m\tau),$ $\operatorname{Im}\tau>0\},$ $E(s)$ is the zeta
function of the Laplacian of the flat metric on the elliptic curves $E_{\tau
},$ the regularized determinant of the Laplacian is $\exp(-\frac{d}%
{ds}E(s)|_{s=0})$ and $\eta^{24}$ is equal to the discriminant of the elliptic
curve $E_{\tau}. $ $\eta^{24}$ vanishes at $\infty,$ which corresponds to an
elliptic curve with the node. Thus the Kronecker limit formula is an explicit
formula for the determinant of the Laplacian of an elliptic curve and gives a
relation between the spectrum of the Laplacian and the discriminant of
elliptic curves. The Kronecker limit formula has a modern interpretation as
the Quillen norm of a section of the determinant line bundle.

There is a simple non formal explanation of the above mentioned fact. It is a
well known fact that the spectrum of the Laplacian of a Riemannian metric on a
compact manifold is discrete. When the manifold acquires singularities than
the specter becomes continuous. This phenomenon suggests that when the metric
''degenerates'' together with the manifold, then the regularized determinant
vanishes on the points that parametrize the singular varieties. The problem is
how to relate the specter of the Laplacian with the discriminant locus. The
relation is suggested by the theory of determinant line bundles on the moduli
space, their Quillen metrics and the Ray-Singer torsion as developed recently
by Quillen, Donaldson, Bismut, Gillet and Soul\'{e} and others.

The problem that we are going to study in a series of two papers is to find
the generalization of the analogue of the Dedekind eta function for odd
dimensional CY manifolds.

The idea on which these two papers are based is very simple. The Quillen
metric is related to the spectral properties of the Laplacian acting on (0,q)
forms in case of K\"{a}hler manifolds. The main question is when the Ray
Singer analytic torsion is the Quillen metric of some holomorphic section of
the determinant line bundle. It is easy to prove that if the index of the
$\overline{\partial}$ operator is zero, then one can construct a non vanishing
$C^{\infty}$ section $\det(\overline{\partial})$ of the determinant line
bundle $\mathcal{L}$ up to a constant whose Quillen norm is exactly the
analytic Ray Singer torsion. We will show that knowing the existence of the
non vanishing section $\det(\overline{\partial})$ implies that there exists a
holomorphic section $\eta^{N}$ of some power of the determinant line bundle
which vanishes on $\mathcal{D}_{\infty}=\overline{\mathcal{M}(M)}$
$\backslash$%
$\mathcal{M}(M),$ where $\overline{\mathcal{M}(M)}$ is some projective
compactification of $\mathcal{M}(M)$ such that $\mathcal{D}_{\infty}%
=\overline{\mathcal{M}(M)}$
$\backslash$%
$\mathcal{M}(M)$ is a divisor with normal crossings. According to Viewheg
$\mathcal{M}(M)$ is a quasi projective variety. See \cite{V}.

For this program we need the analogue of the variational formulas for the
determinant of the Laplacian of a CY metric acting on (0,q) forms. The
variational formulas are very important in the construction of the holomorphic
section $\eta^{N}$ mentioned before.

In this paper we will generalize our variational formulas that were proved for
the Laplacian of a CY metric acting on functions to (0,q) forms. See
\cite{JT95} and\ \cite{JT98}.

We discussed the problem of finding the relations between the spectral
properties of the Laplacian of CY metric on K3 surfaces in a series of joint
papers with J. Jorgenson. (See \cite{JT95}, \cite{JT96}, \cite{JT098} and
\cite{JT98}.)\footnote{The mistakes that appeared in these papers are
corrected in \cite{JT99}.} The results of these papers showed that the problem
of relating the spectral properties of the Laplacian of CY metric on even
dimensional CY manifold is very delicate one. For example, in the case of
algebraic polarized K3 surfaces we showed that the determinant of the
Laplacian defines the discriminant locus of polarized K3 surfaces for
polarization classes $e$ such that the Baily Borel compactification of the
moduli space of algebraic pseudo-polarized K3 surfaces $\Gamma_{n}\backslash
SO_{0}(2,19)/SO(2)\times SO(19)$ contains only one zero dimensional cusp$.$ In
the other cases the discriminant locus can not be recovered from the spectral
properties of the Laplacians of CY metrics. The difficulties in the even
dimensional case are based on the fact that Ray Singer Analytic torsion is
zero and the index of the $\overline{\partial}$ is equal to 2. So we can not
find a non vanishing canonical section of the determinant line bundle.

The analytic torsion for Enriques surfaces is discussed in \cite{HM1} from
point of view of string theory and in \cite{Y} from mathematical point of
view. Based on these two papers one should consider the Enriques surfaces from
the point of view of the spectral properties of the Laplacian of CY metric as
an odd dimensional CY manifolds.

There are relations between the results our results in the series of the two
papers and the results of in \cite{BCOV}. Some of these relations are
discussed in \cite{TO99} The results and the conjectures stated in thes two
papers are related to the results in \cite{HM1}, \cite{HM2} and \cite{MO}.

This article is organized as follows. In \textbf{Section II} we introduce the
basic definition and some notations. In \textbf{Section III }we review
Kodaira-Spencer-Kuranishi deformation theory of Calabi Yau manifolds following
\cite{To89}. In \textbf{Section IV} we introduce some canonical
identifications of different Hilbert spaces on a CY\ manifold. We prove that
some operators are of trace class and compute their traces. In \textbf{Section
V} we proved that $Tr(\exp(-t\Delta_{q}))=\binom{n}{q}Tr(\exp(-t\Delta_{0})),$
where $\Delta_{q}$ is the Laplacian of the Calabi Yau metric acting on (0,q)
forms$.$ Our proof is based on Bochner type technique. In \textbf{Section VI}
we formulate and prove the main Theorem of this article, namely that
$\log(\det(\Delta_{q}))$ is a potential for the Weil -Petersson metric$.$ In
\textbf{Section VII }we gave some applications of the technique and results
that we used in the previous section. We prove that the coefficients $a_{-k}$
in the short term asymptotic expansion of the trace of the heat kernel of CY
metric for $0\leq k\leq\dim_{\mathbb{C}}$M are constants. In \textbf{Section
VIII }show that Ray Singer analytic torsion of CY metric I(M) is bounded.

\begin{acknowledgement}
The author wants to thank G. Moore for many stimulating conversations about
the topic in this paper. These conversations inspired many of the ideas and
the results in this paper. I want to thank my friend J. Jorgenson for
introducing me to the exciting world of determinants of Laplacians. I am
grateful to S. Donaldson, S.-T. Yau, G. Zuckerman, D. Kazhdan, S. Lang, B.
Lian, J. Li, K. Liu, Y. Eliashberg and R. Donagi for their encouragements,
useful comments and support. I want to thank Sinan Unver for his help. Special
thanks to the National Center for Theoretical Sciences (Taiwan) for their
hospitality during the preparation of this article. I want to express my
special thanks to Prof. Chang and Prof. Wang. I want to thank Yale University
for their hospitality and the opportunity to lecture on part of the material
included in this paper. I want to thank G. Moore, D. Mostow and G. Zuckreman
for making a number of useful comments. I want to thank P. Deligne for his
useful and critical remarks.
\end{acknowledgement}

\section{Some Remarks, Notations and Preliminary Results.}

\subsection{Definition of the Regularized Determinant}

Let (M,g) be an n dimensional Riemannian manifold. Let $\Delta_{q}=dd^{\ast
}+d^{\ast}d$ be the Laplacian acting on the space of q forms on M. It is a
well known fact that the spectrum of the Laplacian $\Delta_{q}$ is positive
and discrete. This means that the non zero eigen values of $\Delta_{q}$ are
$0<\lambda_{1}\leq\lambda_{2}\leq...\leq\lambda_{n}\leq...$ We will define the
zeta function of $\Delta_{q}$ as follows: $\zeta_{q}(s)=\sum_{i=1}^{\infty
}\lambda_{i}^{-s}.$ It is a well known fact that $\zeta_{q}(s)$ is a well
defined analytic function for $\operatorname{Re}(s)\gg C,$ it has a
meromorphic continuation in the complex plane and $0$ is not a pole of
$\zeta_{q}(s).$ Then we define $\det(\Delta_{q})=\exp\left(  -\frac{d}%
{ds}\left(  \zeta_{q}(s)\right)  |_{s=0}\right)  .$

\subsection{Definitions and Notations}

Let M be a n-dimensional K\"{a}hler manifold with a zero canonical class.
Suppose that $H^{k}($M,$\mathcal{O}_{\text{M}})=0$ for 1$\leq k<n.$
\textit{Such manifolds are called Calabi-Yau manifolds. }A pair (M,$L$) will
be called a polarized CY manifold if M is a CY manifold and $L\in H^{2}%
($M,$\mathbb{Z}$)\footnote{Notice that $H^{1,1}($M,$\mathbb{R)=}H^{2}%
($M,$\mathbb{R)}$ since $H^{2}$(M,$\mathcal{O}_{\text{M}})=0$ for CY
manifolds.} is a fixed class such that it represents the imaginary part of a
K\"{a}hler metric on M.

Yau's celebrated theorem asserts the existence of a unique Ricci flat
K\"{a}hler metric g on M such that the cohomology class [Im(g)]=$L$. From now
on we will consider polarized CY manifolds of odd dimension. The polarization
class $L$ determines the CY metric g uniquely. We will denote by
$\bigtriangleup_{q}=\overline{\partial}^{\ast}\circ\overline{\partial
}+\overline{\partial}\circ\overline{\partial}^{\ast}$ the associated
Laplacians that act on smooth $(0,q)$ forms on M for $0\leq q\leq n$.
$\overline{\partial}^{\ast}$ is the adjoint operator of $\overline{\partial}$
with respect to the CY metric g. The determinant of $\ $these operators
$\bigtriangleup_{q},$ defined through zeta function regularization, will be
denoted by det$\left(  \bigtriangleup_{q}\right)  .$

The Hodge decomposition theorem asserts that $\Gamma($M,$\Omega^{0,q}%
)=\operatorname{Im}(\overline{\partial})\oplus\operatorname{Im}(\overline
{\partial}^{\ast})$ for $1\leq q\leq\dim_{\mathbb{C}}M-1.$ The restriction of
$\ \bigtriangleup_{q}$ on $\operatorname{Im}(\overline{\partial})$ will be
denoted by $\bigtriangleup_{q}^{^{\prime}}=\overline{\partial}\circ
\overline{\partial}^{\ast}$, and the restriction of $\ \Delta_{q}$ on
$\operatorname{Im}(\overline{\partial}^{\ast})$ will be denoted by
$\bigtriangleup_{q}^{"}=\overline{\partial}^{\ast}\circ\overline{\partial}.$
Hence we have $Tr(\exp(-t\bigtriangleup_{q})=Tr(\exp(-t\bigtriangleup
_{q}^{^{\prime}})+Tr(\exp(-t\bigtriangleup_{q}^{"}).$ This implies that
$\zeta_{q}(s)=\sum_{k=1}^{\infty}\lambda_{k}^{-s}=\zeta_{q}^{^{\prime}%
}(s)+\zeta_{q}^{"}(s),$ where $\lambda_{k}>0$ are the positive eigen values of
$\bigtriangleup_{q}$ and $\zeta_{q}^{^{\prime}}(s)$ \& $\zeta_{q}^{"}(s)$ are
the zeta functions of $\bigtriangleup_{q}^{^{\prime}}$ and $\bigtriangleup
_{q}^{"}.$ From here and the definition of \ the regularized determinant we
obtain that $\log\det(\bigtriangleup_{q})=\log\det(\bigtriangleup
_{q}^{^{\prime}})+\log\det(\bigtriangleup_{q}^{"}).$ It is a well known fact
that the action of $\bigtriangleup_{q}^{^{"}}$ on $\operatorname{Im}%
\overline{\partial}^{\ast}$ is isospectral to the action of $\bigtriangleup
_{q+1}^{^{\prime}}$ on $\operatorname{Im}\overline{\partial},$ which means
that the spectrum of $\bigtriangleup_{q}^{^{"}}$ is equal to the spectrum of
$\bigtriangleup_{q+1}^{^{\prime}}.$ So we have the equality $\det
(\bigtriangleup_{q}^{"})=\det(\bigtriangleup_{q+1}^{^{\prime}}).$

\begin{notation}
Let f be a map from a set A to a set B and let g be a map from the set B to
the set C, then the compositions of those two maps we will denote by f$\circ$g.
\end{notation}

\section{Kodaira-Spencer-Kuranishi Theory for CY}

\subsection{Basic Definitions}

In \cite{Ti} and \cite{To89} was developed the local deformation theory of CY
manifolds. We will review the results in \cite{Ti} and \cite{To89} in this section.

Let M be an even dimensional C$^{\infty}$ manifold. We will say that M has an
almost complex structure if there exists a section $I\in C^{\infty
}(M,Hom(T^{\ast},T^{\ast})$ such that $I^{2}=-id.$ $T$ is the tangent bundle
\ and $T^{\ast}$ is the cotangent bundle on M. This definition is equivalent
to the following one: Let M be an even dimensional C$^{\infty}$ manifold.
Suppose that there exists a global splitting of the complexified cotangent
bundle $T^{\ast}\otimes\mathbf{C}=\Omega^{1,0}\oplus\Omega^{0,1}$, where
$\Omega^{0,1}=\overline{\Omega^{1,0}}.$ Then we will say that M has an almost
complex structure. We will say that an almost complex structure is an
integrable one, if for each point $x\in$M there exists an open set $U\subset$M
such that we can find local coordinates $z^{1},..,z^{n},$ such that
$dz^{1},..,dz^{n}$ \ are linearly independent in each point $m\in U$ and they
generate $\Omega^{1,0}|_{U}.$

\begin{definition}
\label{belt}Let M be a complex manifold. Let $\phi\in\Gamma($M,$Hom(\Omega
^{1,0},\Omega^{0,1}))$, then we will call $\phi$ a Beltrami differential.
\end{definition}

Since $\Gamma($M,$Hom(\Omega^{1,0},\Omega^{0,1}))\backsimeq\Gamma($%
M,$\Omega^{0,1}\otimes T^{1,0})$, we deduce that locally $\phi$ can be written
as follows: $\phi|_{U}=\sum\phi_{\overline{\alpha}}^{\beta}\overline
{dz}^{\alpha}\otimes\frac{\partial}{\partial z^{\beta}}$. From now on we will
denote by

\begin{center}
$A_{\phi}=\left(
\begin{array}
[c]{cc}%
id & \phi(\tau)\\
\overline{\phi(\tau)} & id
\end{array}
\right)  .$
\end{center}

We will consider only those Beltrami differentials $\phi$ such that
det($A_{\phi})\neq0.$ The Beltrami differential \ $\phi$\ defines an
integrable complex structure on M if and only if the following equation holds:
$\overline{\partial}\phi+\frac{1}{2}\left[  \phi,\phi\right]  =0,$ where

\begin{center}
$\left[  \phi,\phi\right]  |_{U}:=\sum_{\nu=1}^{n}\sum_{1\leqq\alpha
<\beta\leqq n}\left(  \sum_{\mu=1}^{n}\left(  \phi_{\overline{\alpha}}^{\mu
}\left(  \partial_{\mu}\phi_{\overline{\beta}}^{\nu}\right)  -\phi
_{\overline{\beta}}^{\mu}\left(  \partial_{\nu}\phi_{\overline{\alpha}}^{\nu
}\right)  \right)  \right)  \overline{dz}^{\alpha}\wedge\overline{dz}^{\beta
}\otimes\frac{\partial}{dz^{\nu}}.$ (See \cite{KM}.)
\end{center}

\subsection{Kuranishi Space and Flat Local Coordinates}

Kuranishi proved the following Theorem:

\begin{theorem}
\label{Kur}Let $\left\{  \phi_{i}\right\}  $ be a basis of harmonic $(0,1)$
forms of $\mathbb{H}^{1}($M$,T^{1,0})$ on a Hermitian manifold M. Let $G$ be
the Green operator and let $\phi(\tau^{1},..,\tau^{N})$ be defined as follows:\ 
\end{theorem}

\begin{center}
$\phi(\tau^{1},..,\tau^{N})=\sum_{i=1}^{N}\phi_{i}\tau^{i}+\frac{1}%
{2}\overline{\partial}^{\ast}G[\phi(\tau^{1},..,\tau^{N}),\phi(\tau
^{1},..,\tau^{N})]$,
\end{center}

\textit{then there exists }$\varepsilon>0$\textit{\ such that if }$\tau
=(\tau^{1},..,\tau^{N})$\textit{\ satisfies }$|\tau_{i}|<\varepsilon$\textit{,
then }$\phi(\tau^{1},..,\tau^{N})$\textit{\ is a\ global }$C^{\infty}%
$\textit{\ section of the bundle }$\Omega^{(0,1)}\otimes T^{1,0}$.(See
\cite{KM}.)

Based on the Theorem \ref{Kur} we proved in \cite{To89} the following Theorem:

\begin{theorem}
\label{tod1}Let M be a CY manifold and let $\left\{  \phi_{i}\right\}  $ be a
basis of harmonic $(0,1)$ forms with coefficients in $T^{1,0}$ of
$\mathbb{H}^{1}($M$,T^{1,0}),$ then the equation: $\overline{\partial}%
\phi+\frac{1}{2}\left[  \phi,\phi\right]  =0$ has a solution in the form:
\end{theorem}

\begin{center}
$\phi(\tau_{1},..,\tau_{N})=\sum_{i=1}^{N}\phi_{i}\tau^{i}+\sum_{|I_{N}%
|\geqq2}\phi_{I_{N}}\tau^{I_{N}}=\sum_{i=1}^{N}\phi_{i}\tau^{i}+\frac{1}%
{2}\overline{\partial}^{\ast}G[\phi(\tau^{1},..,\tau^{N}),\phi(\tau
^{1},..,\tau^{N})],$
\end{center}

\textit{\ }$\overline{\partial}^{\ast}\phi(\tau_{1},..,\tau_{N})=0$,
\textit{where } $I_{N}=(i_{1},..,i_{N})$\ \ \textit{is a multi-index},
\ $\phi_{I_{N}}\in C^{\infty}($M$,\Omega^{0,1}\otimes T^{1,0})$, $\tau^{I_{N}%
}=(\tau^{i})^{i_{1}}..(\tau^{N})^{i_{N}}$ \textit{and for some} $\varepsilon
>0$ \textit{\ } $\phi(\tau)\in C^{\infty}($M,$\Omega^{0,1}\otimes T^{1,0})$
\textit{if} $\ |\tau^{i}|<\varepsilon$ \ \textit{and} $i=1,..,N.$ See
\cite{Ti} and $\cite{To89}.$

It is a standard fact from Kodaira-Spencer-Kuranishi deformation theory that
for each $\tau=(\tau^{1},..,\tau^{N})$ as in Theorem \ref{tod1} the Beltrami
differential $\phi(\tau^{1},..,\tau^{N})$ defines a new integrable complex
structure on M, i.e. the points of $\mathcal{K},$ where $\mathcal{K}%
:\{\tau=(\tau^{1},..,\tau^{N})|$ $|\tau^{i}|<\varepsilon\}$ defines a family
of operators $\overline{\partial}_{\tau}$ on the $C^{\infty}$ family
$\mathcal{K}\times M\rightarrow M,$ parametrized by $\mathcal{K}$ \ and
$\overline{\partial}_{\tau}$ are integrable in the sense of
Newlander-Nirenberg. Moreover it was proved by Kodaira, Spencer and Kuranishi
that we get a complex analytic family of CY manifolds $\pi
:\mathcal{X\rightarrow K},$ where as $C^{\infty}$ manifold
$\mathcal{X\backsimeq K}\times M.$ The family $\pi:\mathcal{X\rightarrow K}$
is called the Kuranishi family. The operators $\overline{\partial}_{\tau}$ are
defined as follows:

\begin{definition}
\label{tod3}Let $\{\mathcal{U}_{i}\}$ be an open covering of M, with local
coordinate system in $\mathcal{U}_{i}$ given by $\{z_{i}^{k}\}$ with
$k=1,...,n=$dim$_{\mathbb{C}}$M. Assume that: $\phi(\tau^{1},..,\tau
^{N})|_{\mathcal{U}_{i}}$ is given by:
\end{definition}

\begin{center}
$\phi(\tau^{1},..,\tau^{N})=\sum_{j,k=1}^{n}(\phi(\tau^{1},..,\tau
^{N}))_{\overline{j}}^{k}d\overline{z}^{j}\otimes\frac{\partial}{\partial
z^{k}}.$
\end{center}

\textit{Then we define \ }$(\overline{\partial})_{\tau,\overline{j}}%
=\frac{\overline{\partial}}{\overline{\partial z^{j}}}-\sum_{k=1}^{n}%
(\phi(\tau^{1},..,\tau^{N}))_{\overline{j}}^{k}\frac{\partial}{\partial z^{k}}.$

\begin{definition}
\label{flat}The coordinates $\tau=(\tau^{1},..,\tau^{N})$ defined in Theorem
\ref{tod1} will be fixed from now on and will be called the flat coordinate
system in $\mathcal{K}$.
\end{definition}

\subsection{Weil-Petersson Metric}

It is a well known fact from Kodaira-Spencer-Kuranishi theory that the tangent
space $T_{\tau,\mathcal{K}\text{ }}$at a point $\tau\in\mathcal{K} $ can be
identified with the space of harmonic (0,1) forms with values in the
holomorphic vector fields $\mathbb{H}^{1}($M$_{\tau},T$). We will view each
element $\phi\in\mathbb{H}^{1}($M$_{\tau},T$) as a pointwise linear map from
$\Omega_{\text{M}_{\tau}}^{(1,0)}$ to $\Omega_{\text{M}_{\tau}}^{(0,1)}.$
Given $\phi_{1}$ and $\phi_{2}\in\mathbb{H}^{1}($M$_{\tau},T$)$,$ the trace of
the map: $\phi_{1}\overline{\phi_{2}}:$ $\Omega_{\text{M}_{\tau}}%
^{(0,1)}\rightarrow\Omega_{\text{M}_{\tau}}^{(0,1)}$ at the point $m\in
$M$_{\tau}$ with respect to the metric g is simply:

\begin{center}
$Tr(\phi_{1}\overline{\phi_{2}})(m)=\sum_{k,l,m,p=1}^{n}(\phi_{1}%
)_{\overline{l}}^{k}(\overline{\phi)_{\overline{p}}^{m}}g^{\overline{l}%
,p}g_{k,\overline{m}}.$
\end{center}

\begin{definition}
\label{WP}We will define the Weil-Petersson metric on $\mathcal{K}$ via the
scalar product:
\end{definition}

\begin{center}
$<\phi_{1},\phi_{2}>=\int_{\text{M}}Tr(\phi_{1}\overline{\phi_{2}})vol(g).$
\end{center}

We proved in \cite{To89} that the coordinates $\tau=(\tau^{1},..,\tau^{N}) $
as defined in Definition \ref{flat} are flat in the sense that the
Weil-Petersson metric is K\"{a}hler and in these coordinates we have that the
components $g_{i,\overline{j}}$ of the Weil Petersson metric are given by the
following formulas in these coordinates:

\begin{center}
$g_{i,\overline{j}}=\delta_{i,\overline{j}}+R_{i,\overline{j},l,\overline{k}%
}\tau^{l}\overline{\tau^{k}}+O(\tau^{3}).$
\end{center}

On page 332 of \cite{To89} the following results is proved:

\begin{lemma}
\label{sym}Let $\phi\in\mathbb{H}^{1}($M$,T$) be a harmonic form with respect
to the CY metric g. Let
\end{lemma}

\begin{center}
$\phi|_{U}=\sum_{k,l=1}^{n}\phi_{\overline{k}}^{l}\overline{dz}^{k}%
\otimes\frac{\partial}{\partial z^{l}},$
\end{center}

\textit{then }$\phi_{\overline{k},\overline{l}}=\sum_{j=1}^{n}g_{j,\overline
{k}}\phi_{\overline{l}}^{j}=\sum_{j=1}^{n}g_{j,\overline{l}}\phi_{\overline
{k}}^{j}=\phi_{\overline{l},\overline{k}}.$

We will use Lemma \ref{sym} to prove the following theorem:

\subsection{Infinitesimal Deformation of the Imaginary Part of the
Weil-Petersson Metric}

\begin{theorem}
\label{const}Near the point $\tau=0$ of the Kuranishi space $\mathcal{K}$ the
imaginary part $\operatorname{Im}(g)$ of the CY metric $g$ has the following
expansion in the coordinates $\tau:=(\tau^{1},..,\tau^{N})$:
$\operatorname{Im}(g)(\tau,\overline{\tau})=\operatorname{Im}(g)(0)+O(\tau^{2}).$
\end{theorem}

\textbf{PROOF: }In \cite{To89} we proved that the forms $\theta_{\tau}%
^{k}=dz^{k}+\sum_{l=1}\phi(\tau^{1},..,\tau^{N})_{\overline{l}}^{k}%
d\overline{z^{l}}$ ($k=1,.,n)$ form a basis of $(1,0)$ forms relative to the
complex structure defined by $\tau\in\mathcal{K}$ in $\mathcal{U\subset}$M. Let

\begin{center}
$\operatorname{Im}(g_{\tau})=\sqrt{-1}\sum_{1\leq k\leq l\leq n}%
g_{k,\overline{l}}(\tau,\overline{\tau})$ $\theta_{\tau}^{k}\wedge
\overline{\theta_{\tau}^{l}}.$
\end{center}

and

\begin{center}
$g_{k,\overline{l}}(\tau,\overline{\tau})=g_{k,\overline{l}}(0)+\sum_{i=1}%
^{N}\left(  \left(  g_{k,\overline{l}}(1)\right)  _{i}\tau^{i}+\left(
g_{k,\overline{l}}^{^{\prime}}(1)\right)  _{i}\overline{\tau^{i}}\right)
+O(2)$.
\end{center}

Substituting in the expression for $\operatorname{Im}(g_{\tau})$ the
expressions for $\theta_{\tau}^{k}$ we get the following formula\textit{:}

\begin{center}
$\operatorname{Im}(g_{\tau})=\sqrt{-1}\sum_{1\leq k\leq l\leq n}%
g_{k,\overline{l}}(\tau,\overline{\tau})\theta_{\tau}^{k}\wedge\overline
{\theta_{\tau}^{l}}=\sqrt{-1}\sum_{1\leq k\leq l\leq n}g_{k,\overline{l}%
}(0)dz^{k}\wedge\overline{dz^{l}}+$

$+\sum_{i=1}^{N}\tau^{i}\sqrt{-1}\left(  \sum_{1\leq k\leq l\leq n}\left(
\left(  g_{k,\overline{l}}(1)\right)  _{i}dz^{k}\wedge\overline{dz^{l}}%
+\sum_{m=1}^{n}(g_{k,\overline{m}}\overline{\phi_{i,\overline{l}}^{m}%
}-g_{l,\overline{m}}\overline{\phi_{i,\overline{k}}^{m}})dz^{k}\wedge
dz^{l}\right)  \right)  $

$+\sum_{i=1}^{N}\overline{\tau^{i}}\overline{\sqrt{-1}\left(  \sum_{1\leq
k\leq l\leq n}\left(  \left(  g_{k,\overline{l}}(1)\right)  _{i}dz^{k}%
\wedge\overline{dz^{l}}+(\sum_{m=1}^{n}(g_{k,\overline{m}}\overline
{\phi_{i,\overline{l}}^{m}}-g_{l,\overline{m}}\overline{\phi_{i,\overline{k}%
}^{m}})dz^{k}\wedge dz^{l}\right)  \right)  }.$
\end{center}

From Lemma \ref{sym} we conclude that $\sum_{m=1}^{n}(g_{k,\overline{m}%
}\overline{\phi_{i,\overline{l}}^{m}}-g_{l,\overline{m}}\overline
{\phi_{i,\overline{k}}^{m}})=0$ and so:

\begin{center}
$\operatorname{Im}(g_{\tau})=\sqrt{-1}\sum_{1\leq k\leq l\leq n}%
g_{k,\overline{l}}(0)dz^{k}\wedge\overline{dz^{l}}+$

$\sum_{i=1}^{N}\tau^{i}\sqrt{-1}\left(  \sum_{1\leq k\leq l\leq n}\left(
g_{k,\overline{l}}(1)\right)  _{i}dz^{k}\wedge\overline{dz^{l}}\right)
+\sum_{i=1}^{N}\overline{\tau^{i}}\sqrt{-1}\overline{\sum_{1\leq k\leq l\leq
n}\left(  g_{k,\overline{l}}(1)\right)  _{i}dz^{k}\wedge\overline{dz^{l}}}+O(2).$
\end{center}

Let us define (1,1) forms\textit{\ }$\psi_{i}:$

\begin{center}
$\psi_{i}=\sqrt{-1}\left(  \sum_{1\leq k\leq l\leq n}\left(  g_{k,\overline
{l}}(1)\right)  _{i}dz^{k}\wedge\overline{dz^{l}}\right)  .$
\end{center}

Since

\begin{center}
[$\operatorname{Im}(g_{\tau})$]=[$\operatorname{Im}(g_{0})+\sum_{i=1}^{N}%
\tau^{i}\psi_{i}+\sum_{i=1}^{N}\overline{\tau^{i}\psi_{i}}+O(\tau
^{2})]=[\operatorname{Im}(g_{0})]$
\end{center}

we deduce that\textit{\ }each\textit{\ }$\psi_{i}$\textit{\ }is an exact
form\textit{, i.e.: }$\psi_{i}=\sqrt{-1}\partial\overline{\partial}f_{i},$
where \textit{\ }$f_{i}$ are globally defined functions on\textit{\ }%
M\textit{. }If we prove that \textit{\ }$\psi_{i}=0$ our theorem will
follow.\textit{\ }In \textit{\cite{To89}\ }we proved that\textit{:
}det$(g_{\tau})= $det$(g_{0})+O(2).$ From this result we deduce by direct
computations that:

\begin{center}
det$(g_{\tau})=$det$(g_{0})+\sum_{i=1}^{N}\tau^{i}\left(  \sqrt{-1}\sum
_{k,l}g^{\overline{l},k}\partial_{k}\overline{\partial_{l}}(f_{i})\right)
+\sum_{i=1}^{N}\overline{\tau^{i}}($complex conjugate)+$O(2). $
\end{center}

Hence we obtain that for each i we have\textit{: }$\sum_{k,l}g^{\overline
{l},k}\partial_{k}\overline{\partial_{l}}(f_{i})=\triangle(f_{i})=0,$ where
$\triangle$ is the Laplacian of the metric g. From the maximum principle, we
deduce that all\textit{\ }$f_{i}$ are constants. Theorem \ref{const} is
proved\textit{. }$\blacksquare.$

\section{Hilbert Spaces and Operators of Trace Class.}

\subsection{Spectral Canonical Identifications of Some Hilbert Spaces}

\begin{definition}
\label{Hilb1}We will denote by $L_{0,q-1}^{2}(\operatorname{Im}(\overline
{\partial}^{\ast}))$ the Hilbert subspace in $L^{2}($M,$\Omega^{(0,q-1)})$
which is the $L^{2}$ completion of \ $\overline{\partial^{\ast}}$ exact forms
in $C^{\infty}($M,$\Omega^{(0,q-1)})$ for $q\geq1.$ In the same manner we will
denote by $L_{0,q}^{2}(\operatorname{Im}(\overline{\partial}))$ the Hilbert
subspace in $L^{2}(\Omega^{(0,q)})$ which is the $L^{2}$ completion of
$\ \overline{\partial}$ exact $(0,q)$ forms in $C^{\infty}($M,$\Omega
^{(0,q)})$ for $q\geq0$ and by $L_{1,q-1}^{2}(\operatorname{Im}(\partial))$ we
denote the Hilbert subspace in $L^{2}(\Omega^{(1,q-1)})$ which is the $L^{2} $
competition of the $\partial$ exact $(1,q-1)$ forms in $C^{\infty}( $%
M,$\Omega^{(1,q-1)})$ . All the completions are with respect to the scalar
product on the bundles $\Omega^{p,q}$ defined by the CY metric g.
\end{definition}

Let $\phi(\tau^{1},..,\tau^{N}$ ) be the solution of the equation
$\overline{\partial}\phi(\tau^{1},..,\tau^{N}$ )$=\frac{1}{2}[\phi(\tau
^{1},..,\tau^{N} $ ),$\phi(\tau^{1},..,\tau^{N}$ $)]$ established in Theorem
\ref{tod1}. From the Definition \ref{belt} of the Beltrami differential we
know that $\phi(\tau^{1},..,\tau^{N}$ ) defines a linear fibrewise map
$\phi(\tau^{1},..,\tau^{N}):\Omega^{(1,0)}\rightarrow\Omega^{(0,1)}$. So

\begin{center}
$\phi(\tau^{1},..,\tau^{N}$ )$\in C^{\infty}(M,Hom(\Omega^{(1,0)}%
,\Omega^{(0,1)}.$
\end{center}

We define the following linear map between the vector bundles $\phi\wedge id:
$ $\Omega^{(1,q-1)}\rightarrow\Omega^{(0,q)}$ as $\phi(dz^{i}\wedge
\alpha)=\phi(dz^{i})\wedge\alpha.$

\begin{definition}
\label{Hilb0} \textit{For each }$1\leq q\leq n,$ $\phi\wedge id$
\textit{\ defines a natural operator }F(q,$\phi)$ \textit{between the Hilbert
spaces} $L^{2}($M, $\Omega^{(1,q-1)})$ and $L^{2}$(M,$\Omega^{(0,q)}$).
\end{definition}

\begin{definition}
\label{Hilb}The restriction of \ the map F(q,$\phi$) on the subspace
$\operatorname{Im}(\partial)\subset L^{2}($M, $\Omega^{(1,q-1)}))$ to
$\operatorname{Im}(\overline{\partial)}\subset L^{2}$(M,$\Omega^{(0,q)}$) will
be denoted by $F^{\prime}(q,\phi)$.
\end{definition}

\begin{lemma}
\label{Hilb2} The Hilbert subspaces $L_{0,q-1}^{2}(\operatorname{Im}%
(\overline{\partial}^{\ast})),$ $L_{0,q}^{2}(\operatorname{Im}(\overline
{\partial}))$ and $L_{1,q-1}^{2}(\operatorname{Im}(\partial))$ are invariant
with respect to the Laplacians
\end{lemma}

\begin{center}
$\triangle_{q-1}^{"}=\overline{\partial}_{q-1}^{\ast}\circ\overline{\partial
}_{q},$ $\triangle_{q}^{^{\prime}}=\overline{\partial}_{q}\circ\overline
{\partial}_{q+1}^{\ast},$
\end{center}

\textit{i.e}.

\begin{center}
$\triangle_{q-1}^{"}(L_{0,q-1}^{2}(\operatorname{Im}(\overline{\partial}%
^{\ast}))=L_{0,q-1}^{2}(\operatorname{Im}(\overline{\partial}^{\ast})),$
$\triangle_{q}^{^{\prime}}(L_{0,q}^{2}(\operatorname{Im}(\overline{\partial
}))=L_{0,q}^{2}(\operatorname{Im}(\overline{\partial}))$ \textit{and
}$\triangle_{q}^{"}(L_{1,q-1}^{2}(\operatorname{Im}(\partial))=L_{0,q-1}%
^{2}(\operatorname{Im}(\partial))$
\end{center}

\textbf{PROOF: }The proof of this lemma is standard fact from K\"{a}hler
geometry. The first two identities followed from directly from the definition
of the Laplacian. The last equality follows from the fact that in K\"{a}hler
geometry the Laplacians $\triangle_{q-1}=\overline{\partial}_{q-1}%
\overline{\partial}_{q}^{\ast}+\overline{\partial}_{q-1}^{\ast}\overline
{\partial}_{q}$ and $\triangle_{q-1}^{^{\prime}}=\partial_{q-1}\partial
_{q}^{\ast}+\partial_{q-1}^{\ast}\partial_{q}$ coincide, $\overline{\partial
}_{q}^{\ast}=[\Lambda,\partial_{q}]$ and $\partial_{q}^{\ast}=[\Lambda
,\overline{\partial}_{q}]$. See \cite{KM} and \cite{RS}. Our lemma is proved.
$\blacksquare.$

Let us denote by $\{\omega_{i}(0,q-1)|$ $i=1,...,\infty\}$ all the eigen forms
of the Laplacian $\triangle_{q-1}$in the Hilbert space $L_{0,q-1}%
^{2}(\operatorname{Im}(\overline{\partial}^{\ast}))$ with norm equal to one.
In the same way we will denote by $\{\omega_{i}(0,q)|$ $i=1,...,\infty\}$ all
the eigen forms of the Laplacian $\triangle_{q}$ with norm one in the Hilbert
space $L_{0,q}^{2}(\operatorname{Im}(\overline{\partial}))$ and by
$\{\omega_{i}(1,q-1)|$ $i=1,...,\infty\}$ all the eigen forms of norm one of
the Laplacian $\triangle_{q}$in the Hilbert space $L_{1,q-1}^{2}%
(\operatorname{Im}(\partial)).$

\begin{lemma}
\label{orth}The forms $\{\omega_{i}(0,q-1)|$ $i=1,...,\infty\},$ $\{\omega
_{i}(0,q)|$ $i=1,...,\infty\}$ and $\{\omega_{i}(1,q-1)|$ $i=1,...,\infty\}$
form orthonormal bases in the Hilbert spaces $L_{0,q-1}^{2}(\operatorname{Im}%
(\overline{\partial}^{\ast})),$ $L_{0,q}^{2}(\operatorname{Im}(\overline
{\partial}))$ and $L_{0,q-1}^{2}(\operatorname{Im}(\partial)).$
\end{lemma}

\textbf{PROOF: }The proof of this lemma is standard fact from the theory of
self-adjoint compact operators in Hilbert spaces. See \cite{Gil}.
$\blacksquare.$

\begin{lemma}
\label{Hilb3}Let (M,g) be a K\"{a}hler manifold with a K\"{a}hler metric g. Let
\end{lemma}

\begin{center}
($L_{0,q-1}^{2}(\operatorname{Im}(\overline{\partial^{\ast}})),\{\omega
_{i}(0,q-1)\})$, $(L_{0,q}^{2}(\operatorname{Im}(\overline{\partial
})),\{\omega_{i}(0,q)\}$\ 
\end{center}

\textit{and }$(L_{1,q-1}^{2}(\operatorname{Im}(\partial)),\{\omega
_{i}(1,q-1)\})$ \textit{be the Hilbert spaces with orthonormal bases defined
in Definition} \ref{Hilb1} \textit{for q}$\geq1$. \textit{Then}

\begin{center}
$\overline{\partial}\left(  \frac{\omega_{i}(0,q-1)}{\parallel\overline
{\partial}\omega_{i}(0,q-1)\Vert^{2}}\right)  =\lambda_{i}\omega_{i}(0,q)$ and
$\partial\left(  \frac{\omega_{i}(0,q-1)}{\parallel\partial\omega
_{i}(0,q-1)\Vert^{2}}\right)  =\lambda_{i}\omega_{i}(1,q-1).$
\end{center}

\textbf{PROOF: }This a standard fact which can be found in \cite{RS}.
$\blacksquare.$

\begin{remark}
\label{ident} Lemma \ref{Hilb3} gives a natural identification of the Hilbert
spaces $L_{0,q-1}^{2}(\operatorname{Im}(\overline{\partial^{\ast}}))$,
$L_{0,q}^{2}(\operatorname{Im}(\overline{\partial}))$\ and $L_{1,q-1}%
^{2}(\operatorname{Im}(\partial))$ because we can choose natural bases of all
these Hilbert spaces by choosing an orthonormal basis consisting of eigen
forms of the Laplacians. We are using the following orthonormal bases to get
the above identifications:
\end{remark}

\begin{center}
$\{\omega_{i}(0,q\},$ $\{\frac{\partial(\omega_{i}(0,q)}{\Vert\partial
(\omega_{i}(0,q)\Vert^{2}}:=e_{i}\}$ and $\{\frac{\overline{\partial}%
(\omega_{i}(0,q)}{\Vert\overline{\partial}(\omega_{i}(0,q)\Vert^{2}}:=f_{i}\}.$
\end{center}

\subsection{Trace Class Operators in Hilbert Spaces}

\textit{We will define the trace of the operator \ }F$^{^{\prime}}$(q,$\phi
$):$L_{1,q-1}^{2}(\operatorname{Im}(\partial))\rightarrow L_{0,q}%
^{2}(\operatorname{Im}(\overline{\partial}))$ (\textit{if this trace exists)
acting on the identified Hilbert spaces as the usual trace of an operator
acting on a Hilbert space. For example we define the trace of the operator
}F$^{^{\prime}}$(q,$\phi$) \textit{with respect to the orthonormal bases}

\begin{center}
$\{\frac{\partial(\omega_{i}(0,q)}{\Vert\partial(\omega_{i}(0,q)\Vert^{2}%
}=\omega_{i}(1,q):=e_{i}\}$ and $\{\frac{\overline{\partial}(\omega_{i}%
(0,q)}{\Vert\overline{\partial}(\omega_{i}(0,q)\Vert^{2}}=\omega
_{i}(0,q+1):=f_{i}\}.$
\end{center}

\begin{theorem}
\label{tr} Let F'(q,$\phi)$ be defined as in Definition \ref{Hilb}, then F
$^{^{\prime}}$(q,$\phi$) are operators of trace class.
\end{theorem}

\textbf{PROOF: }From the Definition \ref{Hilb} of the operators F$^{^{\prime}%
}(q,\phi)$ we know that they are induced by the fibrewise linear maps
$\phi\wedge id$ :$\Omega^{1,q-1}\rightarrow\Omega^{0,q}.$

Since M is a compact manifold we can choose $N_{1,q-1}$ global $C^{\infty}$
forms $\psi_{i}$ of type (1,q-1) such that they span at each point $y\in M,$
the space $\Omega_{y}^{1,q-1}.$ In the same way we can find $N_{0,q}$ forms
$\sigma_{j}$ of type (0,q) such that they span at each point $y\in M,$ the
space $\Omega_{y}^{1,q-1}.$ Without lost of generality we may assume that both
$\psi_{i}$ and $\sigma_{j}$ are linearly independent vectors in the identified
Hilbert spaces $L_{1,q-1}^{2}(\operatorname{Im}(\partial))$ $\&$ $L_{0,q}%
^{2}(\operatorname{Im}(\overline{\partial})).$

Then the maps F(q,$\phi)$:$L_{1,q-1}^{2}(\operatorname{Im}(\partial
))\rightarrow L_{0,q}^{2}(\operatorname{Im}(\overline{\partial}))$ are given
by $N_{1,q-1}\times N_{0,q}$ matrix. So the maps F(q,$\phi)$ are linear
operators between finite dimensional spaces therefore they are of trace class.
Since F$^{^{\prime}}(q,\phi)$ are the restriction of the trace class operators
F(q-1,$\phi)$, we deduce that F$^{^{\prime}}($q,$\phi)$ are of trace class
too.\textbf{\ }Theorem \ref{tr} is proved. $\blacksquare.$

\begin{corollary}
\label{trace}The operator $\overline{\partial}^{-1}\circ F^{^{\prime}}%
(q,\phi)\circ\partial$ is of trace class.
\end{corollary}

\textbf{PROOF: }We have identified the Hilbert spaces

\begin{center}
$L_{0,q-1}^{2}(\operatorname{Im}(\overline{\partial}^{\ast}))$, $L_{0,q}%
^{2}(\operatorname{Im}(\overline{\partial}))$\ and $L_{1,q-1}^{2}%
(\operatorname{Im}(\partial))$
\end{center}

in Remark \ref{ident}.

The operators $\overline{\partial}^{-1}\circ F^{^{\prime}}(q,\phi
)\circ\partial$ which act on $L_{0,q-1}^{2}(\operatorname{Im}(\overline
{\partial^{\ast}}))$ can be considered as a composition of a differential
operator, operators with a smooth kernel and integral operator by using the
above identification. From Proposition \textbf{2.45} page 96 in the book
\cite{BGV} it follows directly that the operator $\overline{\partial}%
^{-1}\circ F^{^{\prime}}(q,\phi)\circ\partial$ is of trace class. Cor.
\ref{trace} is proved. $\blacksquare.$

\begin{theorem}
\label{holder}For $t>0$ and q$\geq1$ the following equality holds
\end{theorem}

\begin{center}
\textbf{\ }$Tr\left(  \exp\left(  -t(\triangle_{q-1}^{^{"}})\right)
\circ\overline{\partial}^{-1}\circ F^{^{\prime}}(q,\phi)\circ\partial\right)
=Tr\left(  \exp\left(  -t(\triangle_{q}^{^{^{\prime}}}\right)  \circ
F^{^{\prime}}(q,\phi)\right)  =\sum_{i=1}^{\infty}\exp(-t\lambda_{i})a_{ii}.$
\end{center}

\textit{where }$\lambda_{i}$ \textit{are eigen values of} $\triangle
_{q-1}^{^{"}}$ \textit{and we have the following expression for the trace: }

\begin{center}
$Tr\left(  F^{^{\prime}}(q,\phi)\right)  =\sum_{i=1}^{\infty}a_{ii}$
\end{center}

\textit{in the orthonormal bases consisting of eigen vectors of the
corresponding Laplacians as defined in Lemma \ref{orth}}$.$

\textbf{PROOF: }Theorem \ref{tr} and Corollary \ref{trace} imply that the
operators $\overline{\partial}^{-1}\circ F^{^{\prime}}(q,\phi)\circ\partial$
and $F^{^{\prime}}(q,\phi)$ are of trace class. The proof of this theorem is
based on the direct computation of the traces of the operators $\overline
{\partial}^{-1}\circ F^{^{\prime}}(q,\phi)\circ\partial$ and $F^{^{\prime}%
}(q,\phi)$ with respect to the standard bases of orthonormal vectors \ 

\begin{center}
$\left\{  \frac{\partial\omega_{i}(0,q-1)}{\left\|  \partial\omega
_{i}(0,q-1)\right\|  ^{2}}=\omega_{i}(1,q-1)\right\}  $ and $\left\{
\frac{\overline{\partial}\omega_{i}(0,q-1)}{\left\|  \overline{\partial}%
\omega_{i}(0,q-1)\right\|  ^{2}}=\omega_{i}(0,q)\right\}  ,$
\end{center}

where $\Delta_{q-1}\left(  \omega_{i}(0,q-1)\right)  =\lambda_{i}\omega
_{i}(0,q-1),$ $\Delta_{q}\omega_{i}(1,q-1)=\lambda_{i}\omega_{i}(1,q-1)$ and
$\Delta_{q}\omega_{i}(0,q)=\lambda_{i}\omega_{i}(0,q).$(See Lemma
$\ref{Hilb3}$)$.$ Let

\begin{center}
$F^{^{\prime}}(q,\phi)(\omega_{i}(1,q-1))=\sum_{j=1}^{\infty}a_{ij}(\omega_{j}(0,q)).$
\end{center}

\begin{lemma}
\label{trace1}We have the following formula:
\end{lemma}

\begin{center}
Tr ($\overline{\partial}^{-1}\circ F^{^{\prime}}(q,\phi)\circ\partial$%
)$=$Tr($F^{^{\prime}}(q,\phi))=\sum_{i=1}^{\infty}a_{ii}$ and q$\geq1.$
\end{center}

\textbf{PROOF: }The operator $\overline{\partial}^{-1}\circ F^{^{\prime}%
}(q,\phi)\circ\partial$ act on the Hilbert space $L_{0,q-1}^{2}%
(\operatorname{Im}(\overline{\partial^{\ast}}))$ with an orthonormal basis of
non zero eigen vectors of the Laplacian $\Delta_{q-1}$ $\{\omega
_{i}(0,q-1)\}.$ Recall that $\left\|  \partial\omega_{i}(0,q-1)\right\|
=\left\|  \overline{\partial}\omega_{i}(0,q-1)\right\|  =\sqrt{\lambda_{i}}.$
So we have

\begin{center}
$\partial\omega_{i}(0,q-1)=\sqrt{\lambda_{i}}\omega_{i}(1,q-1)$ and
$\overline{\partial}\omega_{i}(0,q-1)=\sqrt{\lambda_{i}}\omega_{i}(0,q).$
\end{center}

From the expression $\ F^{^{\prime}}(q,\phi)(\omega_{i}(1,q-1))=\sum
_{j=1}^{\infty}a_{ij}(\omega_{j}(0,q))$ and above equalities we obtain the
following formula for the matrix of the operator $\overline{\partial}%
^{-1}\circ F^{^{\prime}}(q,\phi)\circ\partial$ in the basis $\{\omega_{i}(0,q-1)\}$

\begin{center}
$\overline{\partial}^{-1}\circ F^{^{\prime}}(q,\phi)\left(  \partial
(\omega_{i}(0,q-1)\right)  =\overline{\partial}^{-1}\circ F^{^{\prime}}%
(q,\phi)\left(  \sqrt{\lambda_{i}}\omega_{i}(1,q-1)\right)  =$

$\sqrt{\lambda_{i}}\left(  \overline{\partial}^{-1}\circ\sum_{j=1}^{\infty
}a_{ij}(\omega_{j}(0,q))\right)  =\left(  \sqrt{\lambda_{i}}\right)
\sum_{j=1}^{\infty}a_{ij}(\overline{\partial}^{-1}\omega_{j}(0,q))=$
\end{center}

Substituting in the last formula the expression $\frac{\overline{\partial
}\omega_{j}(0,q-1)}{\sqrt{\lambda_{j}}}=\omega_{j}(0,q)$ we obtain that

\begin{center}
$\overline{\partial}^{-1}\circ F^{^{\prime}}(q,\phi)\left(  \partial
(\omega_{i}(0,q-1)\right)  =\sqrt{\lambda_{i}}\sum_{j=1}^{\infty}%
a_{ij}(\overline{\partial}^{-1}\frac{\overline{\partial}\omega_{j}%
(0,q-1)}{\sqrt{\lambda_{j}}})=.$

$=\sum_{j=1}^{\infty}\frac{\sqrt{\lambda_{j}}}{\sqrt{\lambda_{i}}}%
a_{ij}(\omega_{j}(0,q-1)).$
\end{center}

So

\begin{center}
$Tr\left(  \overline{\partial}^{-1}\circ F^{^{\prime}}(q,\phi)\circ
\partial\right)  =\sum_{i=1}^{\infty}\left\langle \sum_{j=1}^{\infty}%
\frac{\sqrt{\lambda_{j}}}{\sqrt{\lambda_{i}}}a_{ij}(\omega_{j}(0,q-1)),\omega
_{i}(0,q-1)\right\rangle =$

$\sum_{i=1}^{\infty}a_{ii}\sum_{i=1}^{\infty}\left\langle F^{\prime}%
(q,\phi)(\omega_{i}(0,q-1)),\omega_{i}(0,q-1)\right\rangle .$
\end{center}

On the other hand we know from Theorem \ref{tr} that the operator $F^{\prime
}(q,\phi)$ is of trace class. From the canonical identifications of the
Hilbert spaces $L_{1,q-1}^{2}(\operatorname{Im}\partial)$ and $L_{0,q}%
^{2}(\operatorname{Im}\overline{\partial})$ by the orthonormal eigen forms
with a non zero eigen forms we deduce that

\begin{center}
$Tr(F^{\prime}(q,\phi))=\sum_{i=1}^{\infty}\left\langle \ F^{^{\prime}}%
(q,\phi)(\omega_{i}(1,q-1)),\omega_{i}(0,q)\right\rangle =\sum_{i=1}^{\infty
}\left\langle \sum_{j=1}^{\infty}a_{ij}(\omega_{j}(0,q)),\omega_{i}%
(0,q)\right\rangle $.
\end{center}

Lemma \ref{trace1} is proved. $\blacksquare.$

\textbf{The End of the Proof of Theorem }\ref{holder}: The formulas
$\triangle_{q-1}^{^{^{"}}}\left(  \omega_{j}(0,q-1)\right)  =\lambda_{i}%
\omega_{j}(0,q-1)$ imply that

\begin{center}
$\exp(-t(\triangle_{q-1}^{^{^{"}}})\left(  \omega_{j}(0,q-1)\right)
=\exp(-t\lambda_{j})\omega_{j}(0,q-1).$
\end{center}

From the expression $\overline{\partial}^{-1}\circ F^{^{\prime}}(q,\phi
)\circ\partial(\omega_{i}(0,q-1))=\sum a_{ij}\frac{\sqrt{\lambda_{j}}}%
{\sqrt{\lambda_{i}}}\omega_{j}(0,q-1)$ proved in Lemma \ref{trace1} we deduce:

\begin{center}
$Tr\left(  \exp(-t(\triangle_{q-1}^{^{"}})\circ\overline{\partial}^{-1}\circ
F^{^{\prime}}(q,\phi)\circ\partial\right)  =$

$\sum_{i=1}^{\infty}\left\langle \left(  exp(-t(\triangle_{q-1}^{^{"}}%
)\circ\overline{\partial}^{-1}\circ(F^{^{\prime}}(q,\phi)\right)
(\partial\omega_{i}(0,q-1)),\omega_{i}(0,q-1)\right\rangle =$

$\sum_{i=1}^{\infty}\left\langle \sum_{j=1}^{\infty}exp(-t(\lambda_{j}%
)a_{ij}\frac{\sqrt{\lambda_{j}}}{\sqrt{\lambda_{i}}}\omega_{j}(0,q-1)),\omega
_{i}(0,q-1)\right\rangle =$

$\sum_{i=1}^{\infty}a_{ii}\exp(-t\lambda_{i}).$
\end{center}

So we obtain $Tr$($\exp(-t(\triangle_{q}^{^{"}})\circ\frac{\partial}%
{\partial\tau_{i}}\left(  \overline{\partial}_{\tau}\right)  \left|  _{\tau
=0}\right.  )=\sum_{i=1}^{\infty}a_{ii}\exp(-t\lambda_{i}).$ From the
expression $\triangle_{q}^{^{^{\prime}}}\left(  \overline{\partial}\omega
_{i}(0,q-1)\right)  =\lambda_{i}\left(  \overline{\partial}\omega
_{i}(0,q-1)\right)  $ we obtain that $\exp(-t\triangle_{q}^{^{\prime}}%
)\omega_{i}(1,q-1)=\exp(-\lambda_{i}t)\omega_{i}(1,q-1).$ From the formula
$F^{^{\prime}}(q,\phi)(\omega_{i}(1,q-1)))=\sum_{j=1}^{\infty}a_{ij}%
(\omega_{j}(0,q))$ we conclude:

\begin{center}
$\left(  \exp(-t\triangle_{q}^{^{\prime}})\circ F^{^{\prime}}(q,\phi)\right)
(\omega_{i}(1,q-1))=\sum_{j=1}^{\infty}a_{ij}\exp(-\lambda_{j}t)\omega_{i}(0,q)$
\end{center}

and so

\begin{center}
$Tr\left(  \exp(-t\triangle_{q}^{^{\prime}})\circ F^{^{\prime}}(q,\phi
)\right)  =\sum_{j=1}^{\infty}a_{ii}\lambda_{i}\exp(-\lambda_{i}t)=Tr\left(
\exp(-t(\triangle_{q-1}^{^{"}})\circ\overline{\partial}^{-1}\circ F^{^{\prime
}}(q,\phi)\circ\partial\right)  .$
\end{center}

Theorem \ref{holder} is proved. $\blacksquare.$

\begin{corollary}
\label{trace4}We have the following formula:
\end{corollary}

\begin{center}
$Tr\left(  \exp(-t(\triangle_{q-1}^{^{"}})\circ\overline{\partial}^{-1}\circ
F^{^{\prime}}(q,\phi)\circ\partial\right)  Tr\left(  \exp(-t\triangle
_{q}^{^{\prime}})\circ F^{^{\prime}}(q,\phi)\right)  =$

$Tr\left(  \exp(-t\triangle_{q}^{^{\prime}})\circ F^{^{\prime}}(q,\phi
)\right)  =\sum_{j=1}^{\infty}a_{ii}\exp(-\lambda_{i}t).$
\end{center}

\begin{corollary}
\label{trace5}a$_{ii}$ tends to zero with $i\rightarrow\infty$ exponentially fast.
\end{corollary}

Repeating the arguments that we used to prove Theorem \ref{holder} we get the
following Theorem:

\begin{theorem}
\label{holder1}For $t>0$ and q$\geq1$ the following equality holds
\end{theorem}

\begin{center}
$Tr\left(  \exp\left(  -t(\triangle_{q-1}^{^{"}})\right)  \circ\partial
^{-1}\circ\overline{F^{^{\prime}}(q,\phi_{j})}\circ F^{^{\prime}}(q_{i}%
,\phi)\circ\partial\right)  =$

$Tr\left(  \exp\left(  -t(\triangle_{q}^{^{^{\prime}}}\right)  \circ
\overline{F^{^{\prime}}(q,\phi_{j})}\circ F^{^{\prime}}(q_{i},\phi)\right)
=\sum_{j=1}^{\infty}c_{ii}\exp(-\lambda_{i}t)$
\end{center}

\textit{where }$Tr\left(  \overline{F^{^{\prime}}(q,\phi_{j})}\circ
F^{^{\prime}}(q_{i},\phi)\right)  =\sum_{j=1}^{\infty}c_{ii}.$

\section{Bochner's Formulas for CY manifolds.}

We will now give explicit expression for Ray-Singer torsion using \ref{var}.
In order to use \ref{var} we need to have information about the relations
between $Tr(\exp(-\triangle_{q}))$ and $Tr(\exp(-\triangle_{0})).$ We will use
Bochner technique to find these relations.

\begin{theorem}
\label{Bochner}$Tr(\exp(-t\triangle_{q}))=\binom{n}{q}Tr(\exp(-t\triangle_{0})).$
\end{theorem}

\textbf{PROOF: }In order to prove Theorem \ref{Bochner}we will use the
following formulas proved in \cite{KM} on page 119: Let M be a K\"{a}hler
manifold and let $\triangle$ be the Laplacian of a K\"{a}hler metric defined
on (p,q) form

\begin{center}
$\phi=\frac{1}{p!q!}\sum\phi_{i_{1},.,i_{p};\overline{j}_{1},.,\overline
{j}_{q}}dz^{i_{1}}\wedge.\wedge dz^{i_{p}}\wedge\overline{dz}^{j_{1}}%
\wedge.\wedge\overline{dz}^{j_{q}},$
\end{center}

then

\begin{center}
$(\triangle\phi)_{i_{1},.,i_{p};\overline{j}_{1},.,\overline{j}_{q}}%
=-\sum_{i,j}g^{\overline{j},i}\nabla_{i}\overline{\nabla}_{j}\phi
_{i_{1},.,i_{p};\overline{j}_{1},.,\overline{j}_{q}}+$

$+\sum_{k}\sum_{l}\sum_{m,n}R^{m}$ $_{i_{k},\overline{j}_{l}}$ $^{\overline
{n}}\phi_{i_{1},.,i_{k-1},m,i_{k},.i_{p}\overline{j}_{1},.,\overline{j}%
_{l-1},\overline{n},\overline{j}_{l+1},..,,\overline{j}_{q}}-$

$-\sum_{k=1}^{n}\sum_{m}R_{\overline{j}_{k}}$ $^{\overline{m}}\phi
_{i_{1},.,i_{p}\overline{j}_{1},.,\overline{j}_{k-1},\overline{m},\overline
{j}_{k+1},..,,\overline{j}_{q}},$
\end{center}

where $R_{i,\overline{j},k,\overline{l}}$ is the curvature of the K\"{a}hler
metric g,\ $\overline{\nabla}_{j}=\overline{\partial}_{j}$ and $\nabla_{i} $
is the covariant derivative in the direction $\frac{\partial}{\partial z^{i}}
$ and

\begin{center}
$R_{\overline{n}}$ $^{\overline{m}}=\sum_{k=1}^{n}g^{\overline{m}%
,k}R_{k,\overline{m}},$
\end{center}

where $R_{k,\overline{m}}$ is the Ricci curvature. If M is a CY manifold and g
is a CY metric, then $R_{k,\overline{m}}=0.$ When $\phi$ is a form of type
$(0,q),$ then from the above mentioned formulas we obtain that:

\begin{center}
$(\triangle\phi)_{\overline{j}_{1},.,\overline{j}_{q}}=-\sum_{n,m}%
g^{\overline{m},n}\nabla_{n}\overline{\nabla}_{m}\phi_{\overline{j}%
_{1},.,\overline{j}_{q}}.$
\end{center}

On page 110 in \cite{KM}the following formula is proved:

\begin{center}
$(\overline{\partial}^{\ast}\phi)_{\overline{j}_{1},.,\overline{j}_{q}%
}=-(-1)^{p}\sum_{m,n}g^{\overline{m},n}\nabla_{n}\phi_{\overline{j}%
_{1},.,\overline{j}_{q}}.$
\end{center}

Using all these formulas we get

\begin{center}
$\left(  \bigtriangleup\phi\right)  _{\overline{j}_{1},.,\overline{j}_{q}%
}=\overline{\partial}^{\ast}\overline{\partial}(\phi_{\overline{j}%
_{1},.,\overline{j}_{q}})=\bigtriangleup_{0}\left(  \phi_{\overline{j}%
_{1},.,\overline{j}_{q}}\right)  .$
\end{center}

From here Theorem \ref{Bochner} follows directly, i.e. $Tr(\exp(-\triangle
_{q}))=\binom{n}{q}Tr(\exp(-\triangle_{0})).$ Our Theorem is proved.
$\blacksquare.$

\section{Variational Formulas.}

Let $<\phi_{i},\phi_{j}>$ be defined as in Definition \ref{WP}, then

\begin{theorem}
\label{var} The following variational formulas hold for CY manifolds of
complex dimension $n\geq2:$
\end{theorem}

\textbf{i. }$\ -(\frac{d^{2}}{\partial\tau_{i}\overline{\partial\tau_{i}}}%
\log(\det(\triangle_{q}^{"}))(0)=\binom{n-1}{q-1}<\phi_{i},\phi_{j}>$ for
1$\leq q\leq n-1.$

\textbf{ii. }$-(\frac{d^{2}}{\partial\tau_{i}\overline{\partial\tau_{i}}}%
(\log(\det(\triangle_{q}^{"}))(0)=<\phi_{i},\phi_{j}>$ for $q=1$ or $n.$

\subsection{Ideas of the Proof}

Let $q\geq1.$ The proof of Theorem \ref{var} is based on the fact that
$\zeta_{q-1}(s)$ is the Mellin transform of $Tr\exp(-t(\triangle_{q-1}^{^{"}%
}),$ i.e. we have

\begin{center}
$\zeta_{q-1}(s)=\frac{1}{\Gamma(s)}\int_{0}^{\infty}(Tr\exp(-t(\triangle
_{q-1}^{^{"}}))t^{s-1}dt.$
\end{center}

The definition of \textit{det}($\triangle_{q-1}^{^{"}})=-\frac{d}{ds}\left(
\zeta_{q-1}(s)\right)  |_{s=0}$, and the power series expansion of zeta
function $\zeta_{q-1}(s)=\zeta_{q-1}(0)+\frac{d}{ds}\left(  \zeta
_{q-1}(0)\right)  s+...$ suggest that in order to compute $\frac{\partial^{2}%
}{\overline{\partial\tau_{i}}\text{ }\partial\tau_{i}}(\det\left(
\log(\triangle_{q-1}^{^{"}})\right)  )|_{\tau=0}$ we need to compute

\begin{center}
$\frac{d}{ds}\left(  \frac{\partial^{2}}{\overline{\partial\tau_{i}}\text{
}\partial\tau_{i}}\frac{1}{\Gamma(s)}\int_{0}^{\infty}Tr\exp(-t(\triangle
_{q-1}^{^{"}})t^{s-1}dt\right)  \left|  _{s=0,\tau=0}\right.  .$
\end{center}

First we will compute $\frac{d}{d\tau_{i}}Tr\exp(-t(\triangle_{q-1}^{^{"}})$
and will prove that

\begin{center}
$\frac{d}{d\tau_{i}}Tr\exp(-t(\triangle_{q-1}^{^{"}})=t\frac{d}{dt}%
Tr\exp\left(  (-t(\triangle_{q-1}^{^{"}})\circ\overline{\partial}^{-1}\circ
F^{\prime}(q,\phi_{i})\circ\partial\right)  =t\sum_{i=1}^{\infty}\exp
(-\lambda_{i}t)\lambda_{i}a_{ii},$
\end{center}

where $Tr(F^{\prime}(q,\phi_{i}))=\sum_{i=1}^{\infty}a_{ii}.$ By integrating
by parts and following closely the arguments from the book \cite{DK} on page
257-260 we will obtain

\begin{center}
$\frac{\partial}{\partial\tau_{i}}\left(  \frac{1}{\Gamma(s)}\int_{0}^{\infty
}Tr\exp(-t(\triangle_{q-1}^{^{"}})t^{s-1}dt\right)  |_{\tau=0}=\frac{s}%
{\Gamma(s)}\int_{0}^{\infty}Tr\left(  \exp(-t(\triangle_{q-1}^{^{"}}%
)\circ\overline{\partial}^{-1}\circ F^{\prime}(q,\phi_{i})\circ\partial
\right)  t^{s-1}dt.$
\end{center}

From the last formula we will obtain that:

\begin{center}
$\frac{d}{ds}\left(  \frac{\partial^{2}}{\overline{\partial\tau_{i}}\text{
}\partial\tau_{i}}\frac{1}{\Gamma(s)}\int_{0}^{\infty}Tr\exp(-t(\triangle
_{q-1}^{^{"}})t^{s-1}dt\right)  \left|  _{s=0,\tau=0}\right.  =$

$-\underset{t\rightarrow0}{\lim}\frac{\overline{\partial}}{\overline
{\partial\tau_{i}}}Tr\left(  \exp(-t(\triangle_{q-1}^{^{"}})\circ
\overline{\partial}^{-1}\circ F^{\prime}(q,\phi_{i})\circ\partial\right)  .$
\end{center}

Direct computation will show that:

\begin{center}
$\underset{t\rightarrow0}{\lim}\frac{\overline{\partial}}{\overline
{\partial\tau_{i}}}Tr\left(  \exp(-t(\triangle_{q-1}^{^{"}})\circ
\overline{\partial}^{-1}\circ F^{\prime}(q,\phi_{i})\circ\partial\right)  =$

$\underset{t\rightarrow0}{\lim}Tr\left(  \exp(-t(\triangle_{q-1}^{^{"}}%
))\circ\partial^{-1}\circ\overline{F^{\prime}(q,\phi_{j})}\circ F^{\prime
}(q,\phi_{i})\circ\partial)\right)  =$

$Tr\left(  \overline{F^{\prime}(q,\phi_{j})}\circ F^{\prime}(q,\phi
_{i})\right)  =\left\langle \phi_{i},\phi_{j}\right\rangle .$
\end{center}

\subsection{Preliminary Results}

\begin{lemma}
\label{L1}The following formulas are true:
\end{lemma}

\begin{center}
$\frac{\partial}{\partial\tau_{i}}\left(  \overline{\partial_{\tau}}\right)
|_{\tau=0}=-$F$^{^{\prime}}(q,\phi_{i})\circ\partial.$
\end{center}

\textbf{PROOF: }From the expression in Definition \ref{tod3} we conclude that

\begin{center}
$\delta_{i}^{^{\prime}}($ $\overline{\partial_{\tau}})=\frac{\partial
}{\partial\tau_{i}}\left(  \overline{\partial}_{\tau}\right)  =\frac{\partial
}{\partial\tau^{i}}(\frac{\overline{\partial}}{\overline{\partial z^{j}}}%
-\sum_{m=1}^{N}(\tau^{m}\sum_{k=1}^{N}(\phi_{m})_{\overline{j}}^{k}%
\frac{\partial}{\partial z^{k}}))+O(\tau^{2})).$
\end{center}

So

\begin{center}
$\frac{\partial}{\partial\tau_{i}}\left(  \overline{\partial_{\tau}}\right)
|_{\tau=0}=-\sum_{k=1}^{N}(\phi_{i})_{\overline{j}}^{k}\frac{\partial
}{\partial z^{k}}.$
\end{center}

Lemma \ref{L1} follows directly from this expression and the Definition
\ref{Hilb} of F$^{^{\prime}}(q,\phi_{i})$. $\blacksquare.$

\begin{lemma}
\label{L2}$\frac{\partial}{\partial\tau_{i}}\left(  \overline{\partial}_{\tau
}^{\ast}\right)  |_{\tau=0}=0.$
\end{lemma}

\textbf{PROOF: }We know from K\"{a}hler geometry that $(\overline
{\partial_{\tau}})^{\ast}=[\Lambda_{\tau},\partial_{\tau}],$ where
$\Lambda_{\tau}$ is the contraction with \textit{(1,1) }vector filed$:$

\begin{center}
$\frac{\sqrt{-1}}{2}\sum_{k,l=1}^{n}$g$_{\tau}^{\overline{k},l}(\theta_{\tau
}^{l})^{\ast}\wedge(\overline{\theta_{\tau}^{k}})^{\ast}.$
\end{center}

on M$_{\tau}$ and $(\theta_{\tau}^{l})^{\ast}$ is (1,0) vector field on
M$_{\tau}$ dual to the (1,0) form $\theta_{\tau}^{i}=dz^{i}+\sum_{j=1}^{N}%
\tau^{j}(\sum_{k=1}^{n}(\phi_{j})_{\overline{k}}^{i}\overline{dz}^{k})).$
Theorem \ref{const}\textit{\ }implies $\frac{\partial}{\partial\tau_{i}%
}(\Lambda_{\tau})|_{\tau=0}=0.$ On the other hand\textit{\ }$\partial_{\tau}%
$\textit{\ }depends antiholomorphically on $\tau$, i.e. it depends
on\textit{\ }$\overline{\tau}=(\overline{\tau_{1}},.,\overline{\tau_{N}}).$ So
we deduce that\textit{:}

\begin{center}
$\frac{\partial}{\partial\tau_{i}}((\overline{\partial_{\tau}})^{\ast}%
)|_{\tau=0}=\left(  [\frac{\partial}{\partial\tau_{i}}(\Lambda_{\tau
}),\partial_{\tau}]+[\Lambda_{\tau},\frac{\partial}{\partial\tau_{i}}%
(\partial_{\tau})]\right)  |_{\tau=0}=0.$
\end{center}

Lemma \ref{L2} is proved. $\blacksquare.$

\subsection{Computation of $\ $holomorphic derivative}

\begin{theorem}
\label{t11}The following formula is true
\end{theorem}

\begin{center}
$\frac{\partial}{\partial\tau_{i}}\left(  Tr(\exp\left(  -t\triangle
_{\tau,q-1}^{"}\right)  \right)  |_{\tau=0}=-tTr\left(  \frac{d}{dt}\left(
\exp\left(  -t\triangle_{q-1}^{"}\right)  \right)  \circ\overline{\partial
}^{-1}\circ F^{\prime}(q,\phi_{i})\circ\partial_{\tau}\right)  |_{\tau=0}.$
\end{center}

\textbf{PROOF: }Direct computations show that:

\begin{center}
$\frac{\partial}{\partial\tau_{i}}\left(  Tr(\exp\left(  -t\triangle
_{\tau,q-1}^{"}\right)  \right)  =-t\exp\left(  \left(  -t\triangle_{\tau
,q-1}^{"}\right)  \circ\frac{d}{d\tau_{i}}\left(  \triangle_{\tau,q-1}%
^{"}\right)  \right)  |_{\tau=0}.$
\end{center}

Lemma \ref{L1} and \ref{L2} imply that

\begin{center}
$\frac{d}{d\tau_{i}}\left(  \triangle_{\tau,q-1}^{"}\right)  |_{\tau
=0}=\left(  \overline{\partial}_{\tau}^{\ast}\circ\frac{d}{d\tau_{i}}%
\overline{\partial}_{\tau}\right)  |_{\tau=0}=-\left(  \overline{\partial
}^{\ast}\circ F^{\prime}(q,\phi_{i})\circ\partial\right)  |_{\tau=0}.$
\end{center}

Since for CY manifolds of complex dimension $\geq3$ the operators
$\overline{\partial}_{\tau}$ give isomorphisms between the spaces of non
constant functions on $M_{\tau}$ $C^{\infty}(M_{\tau})/\mathbb{C}$ and the
space of $C^{\infty}$ $\overline{\partial}_{\tau}$ closed (0,1) forms on
$M_{\tau}.$ So we have that for $\overline{\partial}_{\tau}^{\ast}$ is well
defined on the space of $C^{\infty}$ $\overline{\partial}_{\tau}$ closed (0,1)
forms on $M_{\tau}$ and we have the following formula on $\operatorname{Im}%
\overline{\partial} $:

\begin{center}
$\overline{\partial}_{\tau}^{\ast}=(\triangle_{\tau,q-1}^{"})\circ
\overline{\partial}_{\tau}^{-1}.$
\end{center}

Using all this information we get by direct substitutions that:

\begin{center}
$\frac{\partial}{\partial\tau_{i}}\left(  Tr(\exp\left(  -t\triangle_{q,\tau
}\right)  \right)  |_{\tau=0}=-tTr\left(  \exp\left(  \left(  -t\triangle
_{\tau,q-1}^{"}\right)  \circ\frac{d}{d\tau_{i}}\left(  \overline{\partial
}_{\tau}^{\ast}\circ\overline{\partial}_{\tau}\right)  \right)  \right)
|_{\tau=0}=$

$-tTr\left(  \exp\left(  -t\triangle_{q-1}^{"}\right)  \circ\left(
\overline{\partial}^{\ast}\circ-F^{\prime}(q,\phi_{i})\circ\partial\right)
\right)  =$

$tTr\left(  \exp\left(  -t\triangle_{q-1}^{"}\right)  \circ\triangle_{q-1}%
^{"}\circ\overline{\partial}^{-1}\circ F(q,\phi_{i})\circ\partial\right)
=-tTr\left(  \frac{d}{dt}\exp\left(  -t\triangle_{q-1}^{"}\right)
\circ\overline{\partial}^{-1}\circ F(q,\phi_{i})\circ\partial\right)  .$
\end{center}

Theorem \textbf{\ref{t11}} is proved. $\blacksquare.$

\begin{lemma}
\label{ip}The following formula is true:
\end{lemma}

\begin{center}
$\frac{\partial}{\partial\tau_{i}}\left(  \frac{1}{\Gamma(s)}\int_{0}^{\infty
}Tr\exp(-t(\triangle_{q-1}^{^{"}})t^{s-1}dt\right)  |_{\tau=0}=\frac{s}%
{\Gamma(s)}\int_{0}^{\infty}Tr\left(  \exp(-t(\triangle_{q-1}^{^{"}}%
)\circ\overline{\partial}^{-1}\circ F^{\prime}(q,\phi_{i})\circ\partial
\right)  t^{s-1}dt.$
\end{center}

\textbf{PROOF: }Theorem \ref{t11} imply that we have:

\begin{center}
$\frac{\partial}{\partial\tau_{i}}\left(  \frac{1}{\Gamma(s)}\int_{0}^{\infty
}Tr\exp(-t(\triangle_{q-1}^{^{"}})t^{s-1}dt\right)  |_{\tau=0}=-\frac
{1}{\Gamma(s)}\int_{0}^{\infty}Tr\left(  \frac{d}{dt}\left(  \exp
(-t(\triangle_{q-1}^{^{"}})\right)  \circ\overline{\partial}^{-1}\circ
F^{\prime}(q,\phi_{i})\circ\partial\right)  t^{s}dt.$
\end{center}

Taking into account that Theorem \ref{holder} implies that

\begin{center}
$\underset{t\rightarrow0}{\lim}\left(  Tr\left(  \frac{d}{dt}\left(
\exp(-t(\triangle_{q-1}^{^{"}})\right)  \circ\overline{\partial}^{-1}\circ
F^{\prime}(q,\phi_{i})\circ\partial\right)  \right)  =\underset{t\rightarrow
0}{\lim}\left(  -t\sum_{i=1}^{\infty}\exp(-\lambda_{i}t)\lambda_{i}%
a_{ii}\right)  =0$
\end{center}

and by integrating by parts as in \cite{DK} we derived the formula stated in
Lemma \ref{ip}. Lemma \ref{ip} is proved. $\blacksquare.$

\begin{corollary}
\label{ip1}The following formula is true:
\end{corollary}

\begin{center}
$\frac{d}{ds}\left(  \frac{\partial}{\partial\tau_{i}}(\zeta_{q-1}(s))\right)
|_{s=0}=\underset{t\rightarrow0}{\lim}Tr\left(  \exp(-t(\triangle_{q-1}^{^{"}%
})\circ\overline{\partial}^{-1}\circ F^{\prime}(q,\phi_{i})\circ
\partial\right)  .$
\end{center}

\textbf{PROOF: }Lemma \ref{ip} implies:

\begin{center}
$\frac{d}{d\tau_{i}}\left(  \zeta_{q-1}(s)\right)  =\frac{s}{\Gamma(s)}%
\int_{0}^{\infty}Tr\left(  \exp(-t(\triangle_{q-1}^{^{"}})\circ\overline
{\partial}^{-1}\circ F^{\prime}(q,\phi_{i})\circ\partial\right)  t^{s-1}dt.$
\end{center}

We already computeted the trace of the operator \ $\exp(-t(\triangle
_{q-1}^{^{"}})\circ\overline{\partial}^{-1}\circ F^{\prime}(q,\phi_{i}%
)\circ\partial$ so we obtain:

\begin{center}
$\underset{t\rightarrow0}{\lim}Tr\left(  \exp(-t(\triangle_{q-1}^{^{"}}%
)\circ\overline{\partial}^{-1}\circ F^{\prime}(q,\phi_{i})\circ\partial
\right)  =\underset{t\rightarrow0}{\lim}\sum\exp(-t\lambda_{i})\lambda
_{i}a_{ii}=$

$\sum\lambda_{i}a_{ii}<\infty.$
\end{center}

From the fact that $\frac{s}{\Gamma(s)}=s^{2}+O(s^{3})$ and direct easy
computations we conclude that

\begin{center}
$\frac{d}{ds}\left(  \frac{\partial}{\partial\tau_{i}}(\zeta_{q-1}(s))\right)
|_{s=0}=\underset{t\rightarrow0}{\lim}Tr\left(  \exp(-t(\triangle_{q-1}^{^{"}%
})\circ\overline{\partial}^{-1}\circ F^{\prime}(q,\phi_{i})\circ
\partial\right)  .$
\end{center}

Corollary \ref{ip1} is proved. $\blacksquare.$

\subsection{Computation of the Antiholomorphic Derivative}

Corollary \ref{ip1} implies that we need to compute the antiholomorphic
derivative $\frac{\overline{\partial}}{\overline{\partial\tau}_{j}}Tr\left(
\exp(-t(\triangle_{q-1}^{^{"}})\circ\overline{\partial}^{-1}\circ F^{\prime
}(q,\phi_{i})\circ\partial\right)  $ in order to finish the proof of Theorem
\ref{var}. The computations of the antiholomorphic derivative are based on the
arguments of Quillen as modified in \cite{DK}.

\begin{definition}
\label{ker}We define the function k$_{\tau}^{\#}(w,z,t)$ in a neighborhood of
the diagonal M in MxM as follows: Let $\rho_{\tau}$ be the injectivity radius
on M$_{\tau}.$ Let $d_{\tau}(w,z)$ be the distance between the points $w$ and
$z$ on M$_{\tau}$ with respect to CY metric g$_{\tau}$ and let $\mathcal{P}%
_{\tau}(w,z)(q)$ be the parallel transport of the bundle $\Omega_{\tau
}^{0,q+1}$along the minimal geodesic joining the point $w$ and $z$ with
respect to the Levi Cevita connection of the CY metric$.$ We suppose that
$|\tau|<\varepsilon.$ Let $\delta$ be such that $\delta>\rho_{\tau}.$ Then we
define the function k$_{\tau}^{\#}(w,z,t)$ as a C$^{\infty}$ function using
partition of unity as follows:
\end{definition}

\begin{center}
k$_{\tau}^{\#}(w,z,t)=\left\{
\begin{array}
[c]{ll}%
(4\pi t)^{-\frac{n}{2}}\exp\left(  -\frac{d_{\tau}^{2}(w,z)}{4t}\right)
\mathcal{P}_{\tau}(w,z)(q) & if\text{ }d_{\tau}(w,z)<\rho_{\tau}\\
0 & if\text{ }d_{\tau}(w,z)>\delta.
\end{array}
\right.  $
\end{center}

It was proved in \cite{BGV} on page 87 that we can represent the operator
$\exp(-t\triangle_{\tau,q})$ by an integral kernel k$_{t}(w,z,\tau)$ where

\begin{center}
k$_{t}(w,z,\tau)=(4\pi t)^{-\frac{n}{2}}\exp\left(  -\frac{d_{\tau}^{2}%
(w,z)}{4t}\right)  \left(  \mathcal{P}_{\tau}(w,z)(q)+O(t)\right)  .$
\end{center}

We will denote by M$_{\Delta}\subset$M$\times$M the diagonal in M$\times$M.
Following the arguments from page 258 of \cite{DK}, we will prove the
following theorem:

\begin{theorem}
\label{delta}The following formula holds:
\end{theorem}

\begin{center}
$\frac{\partial}{\partial\tau_{i}}\left(  \log(\det\triangle_{\tau,q}%
^{"}\right)  |_{\tau=0}=\underset{t\rightarrow0}{\lim}\int_{\text{M}}\left(
Tr\left(  \left(  k_{\tau}^{\#}(w,z,t)\left|  _{\operatorname{Im}%
(\overline{\partial})}\right.  \right)  \circ F^{^{\prime}}(q,\phi
_{i})\right)  |_{\tau=0}\right)  vol(g(0)).$
\end{center}

\textbf{PROOF: }An easy calculation, using the fact that $\frac{\overline
{\partial}}{\overline{\partial\tau_{i}}}(F^{^{\prime}}(\phi,q))|_{\tau=0}=0,$
Theorem \ref{holder}, and the definition of $\varepsilon(w,z,t)=\exp
(-\triangle_{q})-$k$_{0}^{\#}(w,z,t),$ show that on the diagonal of MxM for
$t>0$ we have

\begin{center}
$\frac{\partial}{\partial\tau_{i}}\log\left(  \det\triangle_{\tau,q}%
^{"}\right)  |_{\tau=0}=$

$-\underset{t\rightarrow0}{\lim}\left(  \int_{\text{M}_{\Delta}}\left(
Tr\left(  \exp(-t\left(  \triangle_{\tau,q}\right)  \left|
_{\operatorname{Im}(\overline{\partial})}\right.  F^{^{\prime}}(q,\phi
_{i})\right)  |_{\tau=0}\right)  vol\right)  =$

$-\underset{t\rightarrow0}{\lim}\left(  \int_{\text{M}_{\Delta}}\left(
Tr\left(  k_{\tau}^{\#}(w,z,t)\left|  _{\operatorname{Im}(\overline{\partial
})}\right.  \right)  \circ F^{^{\prime}}(q,\phi_{i})\right)  vol\right)  +$

$-\underset{t\rightarrow0}{\lim}\int_{\text{M}}\varepsilon_{0}(w,z,t)\left|
_{\operatorname{Im}(\overline{\partial})}\right.  F^{^{\prime}}(q,\phi_{i})vol=$

$-\underset{t\rightarrow0}{\lim}\left(  (4\pi t)^{-\frac{n}{2}}\int
_{\text{M}_{\Delta}}Tr\left(  \left(  \exp\left(  -\frac{d_{\tau}^{2}%
(w,z)}{4t}\right)  \mathcal{P}_{\tau}(w,z)(q)|_{\tau=0}\right)  \left|
_{\operatorname{Im}(\overline{\partial})}\right.  \circ F^{^{\prime}}%
(q,\phi_{i})\right)  vol\right)  -$

$-\underset{t\rightarrow0}{\lim}\int_{\text{M}_{\Delta}}Tr\left(  \left(
\varepsilon_{0}(w,z,t)\left|  _{\operatorname{Im}(\overline{\partial}%
)}\right.  \right)  \circ F^{^{\prime}}(q,\phi_{i})\right)  vol.$
\end{center}

On the other hand, the definition of k$_{\tau}^{\#}(w,z,t)$ implies that
$\exp(-\triangle_{q})-$k$_{0}^{\#}(w,z,t)=\varepsilon_{0}(w,z,t)$ is bounded
and tends to zero away from the diagonal, as $t$ tends to zero. From here we
deduce that

\begin{center}
$\underset{t\rightarrow0}{\lim}\int_{\text{M}_{\Delta}}Tr\left(  \left(
\varepsilon_{0}(w,z,t)\left|  _{\operatorname{Im}(\overline{\partial}%
)}\right.  \right)  \circ F^{^{\prime}}(q,\phi_{i})\right)  vol=0.$
\end{center}

uniformly in z. Thus, to calculate the limit

\begin{center}
$-\underset{t\rightarrow0}{\lim}\int_{\text{M}_{\Delta}}Tr\left(
\exp(-t\triangle_{\tau,q})\left|  _{\operatorname{Im}(\overline{\partial}%
)}\right.  \circ F^{^{\prime}}(q,\phi_{i})\right)  |_{\tau=0}vol$
\end{center}

we may replace the Heat kernel $\exp(-t\triangle"_{\tau,q})$ by its explicit
approximation k$_{\tau}^{\#}(w,z,t).$ So we deduce that

\begin{center}
$\frac{\partial}{\partial\tau_{i}}\left(  \log(\det\triangle_{\tau,q}%
^{"}\right)  |_{\tau=0}=-\underset{t\rightarrow0}{\lim}\int_{\text{M}_{\Delta
}}Tr\left(  k_{\tau}^{\#}(w,z,t)\left|  _{\operatorname{Im}(\overline
{\partial})}\right.  \circ F^{^{\prime}}(q,\phi_{i})\right)  |_{\tau=0}vol.$
\end{center}

This proves Theorem \ref{delta}. $\blacksquare.$

\begin{corollary}
\label{delta1}Let $\Pr_{q}$ be the projection operator from $L^{2}($%
M,$\Omega^{0,q})$ to $L_{(0,q)}^{2}(\operatorname{Im}(\overline{\partial}))$,
then we have the following formula:
\end{corollary}

\begin{center}
$\frac{d^{2}}{\partial\tau_{i}\overline{\partial\tau_{i}}}\left(  \log
(\det_{q}^{"}(0)\right)  =$

$\underset{t\rightarrow0}{\lim}\frac{1}{\left(  4\pi t\right)  ^{n}}%
\int_{\text{M}}Tr\left(  \left(  \Pr_{q}\left(  \exp-\frac{d_{\tau}^{2}%
(w,z)}{4t}\right)  \circ\delta_{j}^{"}\left(  \mathcal{P}_{\tau}%
(w,z)(q)\left|  _{\tau=0}\right.  \right)  \right)  \circ F^{^{\prime}}%
(q,\phi_{i})\right)  vol.$
\end{center}

\textbf{PROOF: }The proof of the corollary follows directly from Theorem
\ref{delta} and the fact that computation of the trace of a kernel means to
restrict the kernel to the diagonal. From here we deduce that: $\frac
{\overline{\partial}}{\overline{\partial\tau_{j}}}\left(  d_{\tau}%
^{2}(w,z)\right)  |_{w=z}=0.$ Now Corollary\ref{delta1} follows
directly.\textit{\ }$\blacksquare.$

\begin{theorem}
\label{delta2}$\frac{\overline{\partial}}{\overline{\partial\tau_{j}}%
}\mathcal{P}_{\tau}(0,z)(q)=\overline{\text{F}^{\text{'}}\text{(q},\phi_{j}%
)}+O(\overline{\tau}),$where F$^{^{\prime}}$(q,$\phi_{j}$) is defined in
Definition \ref{Hilb}.
\end{theorem}

\textbf{PROOF}: \ We will prove the theorem first for q=1. In this case the
operator $\mathcal{P}_{\tau}$($w,z)$ is the parallel transportation for the
bundle $\Omega^{0,1}$ and it defines a linear map: $\mathcal{P}_{\tau
}(w,a):\Omega_{w,\tau}^{0,1}\rightarrow\Omega_{z,\tau}^{0,1}.$ Once we prove
Theorem \ref{delta2} for q=1, the general case will follow directly from
standard facts from linear algebra.

Since g$_{\tau}$ is a K\"{a}hler metric, the parallel transport operator
$\mathcal{P}_{\tau}(w,z)$ preserves the splitting of the complexified
cotangent bundle of M into (1,0) and (0,1) forms. So the operator
$\mathcal{P}_{\tau}(w,z)$ maps linearly $\Omega_{w,\tau}^{0,1}$ to
$\Omega_{z,\tau}^{0,1}.$ The parallel transportation operators $\mathcal{P}%
_{\tau}(w,z)$ are defined by the Levi-Chevita connection $\nabla_{\tau}$ of
the metrics g$_{\tau}$. We are going to study the local expansion of
$\mathcal{P}_{\tau}(w,z)$ in terms of $\ \tau.$ It is a standard fact that

\begin{center}
$\nabla_{\tau}=\nabla_{\tau}^{1,0}+\nabla_{\tau}^{0,1}=(\partial_{\tau
}-\left(  g_{\tau}^{-1}\partial_{\tau}g_{\tau}\right)  )+\overline{\partial
}_{\tau}.$
\end{center}

In order to define the parallel transportation between $\Omega_{w,\tau}^{0,1}
$ to $\Omega_{z,\tau}^{0,1}$ we need to join the points $w$ and $z$ by
geodesics. We will suppose that $w$ and $z$ are ''close''. This assumption can
be made since we need to compute a trace of an operator given by some kernel.
So our computations will be done on the diagonal M$\subset$M$\times$M. So from
here it follows that we can joint $w$ and $z$ with a unique geodesic. The
parallel transportation of the (0,1) form $\eta\in\Omega_{\tau=0}^{0,1}%
|_{\tau\in\text{N}},$ from a point $w$ to a point $z $ is given by solving the
equations for fix $\tau:$

\begin{center}
$\nabla_{\tau}^{0,1}\left(  \pi_{\tau}^{(0,1)}\eta\right)  :=\overline
{\partial_{\tau}}(\pi_{\tau}^{(0,1)}\eta)=0$

$\nabla_{\tau}^{1,0}\left(  \pi_{\tau}^{(0,1)}\eta\right)  =(\partial_{\tau
}-\left(  g_{\tau}^{-1}\partial_{\tau}g_{\tau}\right)  \left(  \pi_{\tau
}^{(0,1)}\eta(t)\right)  =0$ $\&$ $\eta(0)=\eta.$
\end{center}

where $\pi_{\tau}^{(1,0)}$ and $\pi_{\tau}^{(0,1)}$ are the projection
operators on (1,0) and (0,1) forms on the complex manifold M$_{\tau}$. Without
loss of generality we can assume that $w=0.\,$\ So we can write the following
expression for the parallel transportation operator:

\begin{center}
$\mathcal{P}_{\tau}(0,z)(\eta)=\pi_{\tau}^{(1,0)}(\eta)+\pi_{\tau}%
^{(0,1)}(\eta)+z(\mathcal{B}_{\tau}(\eta))+O(z^{2})$
\end{center}

for point $z\in$M$_{\tau}$ near the fix point $0\in$M$_{\tau}$ and
$\mathcal{B}_{\tau}$ is a linear operator depending on $\tau$, i.e.:
$\mathcal{B}_{\tau}$:$\Omega_{0,\tau}^{(0,1)}\rightarrow\Omega_{z,\tau
}^{(0,1)}.$ Kodaira-Spencer deformation theory implies that $\Omega_{\tau
}^{1,0}$ depends holomorphically on $\tau.$ This fact implies that:
$\delta_{j}^{"}\left(  \pi_{\tau}^{(1,0)}\right)  =0.$ So we obtain the
following formula:

\begin{center}
$\frac{\overline{\partial}}{\overline{\partial\tau_{j}}}(\mathcal{P}_{\tau
}(0,z)(\eta))=\delta_{j}^{"}(\pi_{\tau}^{(0,1)}(\eta))+O(z)$
\end{center}

for $z$ near 0. It is easy to see from the definition of the tangent space to
a point of the Grassmanian of $\Omega_{0,\tau}^{(0,1)}\subset T_{0}^{\ast
}\otimes\mathbb{C}$ that $\delta_{j}^{"}(\pi_{\tau}^{(0,1)})\left|  _{\tau
=0}\right.  =\overline{\phi_{j}}.$ This implies that for $z=0$ we get

\begin{center}
$\frac{\overline{\partial}}{\overline{\partial\tau_{i}}}(\mathcal{P}_{\tau
}(0,z)(\eta)\left|  _{\tau=0}\right.  =\frac{\overline{\partial}}%
{\overline{\partial\tau_{i}}}(\pi_{\tau}^{(0,1)})(\eta)\left|  _{\tau
=0}\right.  =\overline{\phi_{j}}(\eta),$
\end{center}

where $\phi_{j}:C^{\infty}(M,\Omega^{1,0})\rightarrow C^{\infty}%
(M,\Omega^{0,1})$ is the Beltrami operator when $z=0.$ So our Theorem
\ref{delta2} is proved for q=1.

In order to prove Theorem \ref{delta2} for any $q>1,$ we to notice that for
$n_{1}<..<n_{q}:$

\begin{center}
$\nabla\left(  \overline{dz}^{n_{1}}\wedge.\wedge\overline{dz}^{n_{q}}\right)
=\sum_{j=1}^{n}(-1)^{j-1}\left(  \overline{dz}^{n_{1}}\wedge.\wedge\left(
\nabla(\overline{dz}^{n_{j}})\right)  \wedge.\wedge\overline{dz}^{n_{q}%
}\right)  .$
\end{center}

From here the last formula in Theorem \ref{delta2} follows directly once it is
established for the case q=1. Our theorem is proved.$\blacksquare.$

We are now ready to end the proof of Theorem \ref{var}.

\begin{theorem}
\label{var1}The following formula is true for q%
$>$%
0:
\end{theorem}

\begin{center}
$-\frac{d^{2}}{\partial\tau_{i}\overline{\partial\tau_{i}}}\left(  \log
\det\Delta_{\tau,q-1}^{"}\right)  |_{\tau=0}=Tr\left(  F^{^{\prime}}%
(q,\phi_{i})\circ\overline{F^{^{\prime}}(q,\phi_{j})}\right)  .$
\end{center}

\textbf{PROOF}: Corollary \ref{delta1} and Theorem \ref{delta2} implies that
we have the following formula:

\begin{center}
$-\frac{d^{2}}{\partial\tau_{i}\overline{\partial\tau_{i}}}\left(  \log
\det\Delta_{\tau,q-1}^{"}\right)  |_{\tau=0}=$

$-\underset{t\rightarrow0}{\lim}\frac{1}{\left(  4\pi t\right)  ^{n}}%
\int_{\text{M}_{\Delta}}\left(  Tr\left(  \Pr_{q}\left(  \left(  \exp
-\frac{d_{\tau}^{2}(w,z)}{4t}\right)  \circ\delta_{j}^{"}\left(
\mathcal{P}_{\tau}(w,z)(q)\right)  \right)  \left|  _{\tau=0}\right.  \right)
\circ F^{^{\prime}}(q,\phi_{i})\right)  vol=$

$=-\underset{t\rightarrow0}{\lim}\left(  \frac{1}{\left(  4\pi t\right)  ^{n}%
}\int_{\text{M}_{\Delta}}Tr\left(  \Pr_{q}\left(  \exp-\frac{d_{\tau}%
^{2}(w,z)}{4t}\right)  \circ\overline{F^{^{\prime}}(q,\phi_{j})}\circ
F^{^{\prime}}(q,\phi_{i})\right)  vol\right)  .$
\end{center}

Since

\begin{center}
$-\underset{t\rightarrow0}{\lim}\left(  \frac{1}{\left(  4\pi t\right)  ^{n}%
}Tr\left(  \Pr_{q}\left(  \exp-\frac{d_{\tau}^{2}(w,z)}{4t}\right)  \right)
\right)  \left|  _{w=z}\right.  =Dirac(\delta)(z),$
\end{center}

where $Dirac(\delta)$ is the Dirac delta function. So we obtain that

\begin{center}
$-\frac{d^{2}}{\partial\tau_{i}\overline{\partial\tau_{i}}}\left(  \log\left(
\det_{q-1}^{"}(0)\right)  \right)  |_{\tau=0}=Tr\left(  F^{^{\prime}}%
(q,\phi_{i})\circ\overline{F^{^{\prime}}(q,\phi_{j})}\right)  $
\end{center}

This proves Theorem \ref{var1}. $\blacksquare.$

In order to end the proof of the Theorem \ref{var} we will need the following
lemma from linear algebra:

\begin{lemma}
\label{la}Let F be a linear map of a vector space V of dimension n. Then the
linear operator F$\wedge id$ on $\wedge^{q}$V for 1$<q\leq n$ has a trace
given by the formula: Tr(F$\wedge id$)$=\binom{n-1}{q-1}Tr($F).
\end{lemma}

\textbf{PROOF: }The proof is obvious. $\blacksquare.$

\begin{theorem}
$Tr\left(  \text{F}^{^{\prime}}(q,\phi_{i})\overline{\text{F}^{^{\prime}%
}(q,\phi_{j})}\right)  =\binom{n-1}{q-1}<\phi_{i},\phi_{j}>.$
\end{theorem}

\textbf{PROOF: }Applying the variational formula from Theorem \ref{var1} for
q=n-1 we get that

\begin{center}
$-\frac{d^{2}}{\partial\tau_{i}\overline{\partial\tau_{i}}}\left(  \log\left(
\det_{n-1}^{"}(0)\right)  \right)  |_{\tau=0}=Tr\left(  \text{F}^{^{\prime}%
}(n,\phi_{i})\overline{F^{^{\prime}}(n,\phi_{j})}\right)  =Tr\left(
\text{F}^{^{\prime}}(n,\phi_{i})\overline{\text{F}^{^{\prime}}(n,\phi_{j}%
)}\right)  .$
\end{center}

It is easy to see that the composition of the maps F$^{^{\prime}}(n,\phi_{i})
$ and $\overline{F^{^{\prime}}(n,\phi_{j})}$ is defined on the Hilbert space
L$^{2}(M,\Omega^{0,n}),$ i.e.

\begin{center}
F$^{^{\prime}}(n,\phi_{i})\circ\overline{F^{^{\prime}}(n,\phi_{j})}:$
L$^{2}(M,\Omega^{0,n})\rightarrow$L$^{2}(M,\Omega^{0,n}).$
\end{center}

From Hodge theorem and the definition of CY manifold we deduce that:
$L^{2}(M,\Omega^{0,n})=\operatorname{Im}(\overline{\partial})\oplus
\mathbb{C}\overline{\omega_{\tau}}.$ From here and from Lemma \ref{la} we
obtain that

\begin{center}
$-\frac{d^{2}}{\partial\tau_{i}\overline{\partial\tau_{i}}}\left(  \log\left(
\det_{n-1}^{"}(0)\right)  \right)  |_{\tau=0}=Tr\left(  F(1,\phi_{i}%
)\circ\overline{\text{F(1,}\phi_{j})}\right)  .$
\end{center}

On the other hand, the Hodge star operator $\ast$ for CY metric on CY manifold
M, gives us a spectral isomorphism between the Hilbert spaces

\begin{center}
$\ast:L^{2}($M,$\Omega^{0,0})\rightarrow L^{2}($M,$\Omega^{0,n}).$
\end{center}

Since the antiholomorphic form $\overline{\omega_{\tau}}$ is a parallel form
with respect to the Levi-Chevita connection of the CY metric we can deduce
that the Hodge $\ast$ operator gives a spectral isometry between those two
spaces. From this fact and the fact that

\begin{center}
$Tr\left(  F(1,\phi_{i})\circ\overline{\text{F(1,}\phi_{j})}\right)
.=<\phi_{i},\phi_{j}>$
\end{center}

we conclude that:

\begin{center}
$-\frac{d^{2}}{\partial\tau_{i}\overline{\partial\tau_{i}}}\left(  \log\left(
\det_{0}^{"}(0)\right)  \right)  |_{\tau=0}=-\frac{d^{2}}{\partial\tau
_{i}\overline{\partial\tau_{i}}}\left(  \log\left(  \det_{n-1}^{"}(0)\right)
\right)  |_{\tau=0}=Tr\left(  F(1,\phi_{i})\circ\overline{\text{F(1,}\phi
_{j})}\right)  .=<\phi_{i},\phi_{j}>.$
\end{center}

This proves Theorem \ref{var} when q=0 and q=n-1. The formula we just proved, i.e.

\begin{center}
$Tr\left(  F(1,\phi_{i})\circ\overline{\text{F(1,}\phi_{j})}\right)
.=<\phi_{i},\phi_{j}>$
\end{center}

Theorem \ref{var1} and Lemma \ref{la} directly imply Theorem \ref{var} for any
q. $\blacksquare.$

\section{Some Applications.}

\subsection{Computation of the Analytic Torsion}

\begin{theorem}
\label{IZ}Let M be an odd dimensional CY manifold, then
\end{theorem}

\begin{center}
$\log(I(M))=-2\log(\det(\triangle_{0})).$
\end{center}

\textbf{PROOF: }It is a standard fact that

\begin{center}
$\log(I($M$))=\sum_{q=0}^{n=2m+1}(-1)^{q}q\log(\det(\triangle_{q}))=\sum
_{q=1}^{n=2m+1}(-1)^{q}\log(\det(\triangle_{q}^{^{\prime}}))$
\end{center}

From the formulas:

\begin{center}
$\log(\det(\triangle_{q}))=-\zeta_{\triangle_{q}}^{\prime}(0),$

$\zeta_{\triangle_{q}}(s)=\frac{1}{\Gamma(s)}\int_{\text{0}}^{\infty}%
Tr(\exp(-t\triangle_{q}))t^{s-1}dt$

$\zeta_{\triangle_{q}^{^{\prime}}}(s)=\frac{1}{\Gamma(s)}\int_{\text{0}%
}^{\infty}Tr(\exp(-t\triangle_{q}^{^{\prime}}))t^{s-1}dt$
\end{center}

and Theorems \ref{Bochner} \& \ref{var} we deduce that $\log(\det
(\triangle_{q}^{^{\prime}}))=\binom{n-1}{q}\log(\det(\triangle_{0})).$ From
here we obtain that

\begin{center}
$\log(I($M$))=\sum_{q=1}^{n-1}(-1)^{q}\binom{n-1}{q}\log(\det(\triangle
_{0}))+(-1)^{n+1}\log(\det(\triangle_{0})).$
\end{center}

From the equality:

\begin{center}
$(1-1)^{n-1}=\left(  \sum_{q=0}^{n-1}(-1)^{q}\binom{n-1}{q}\right)  =0$
\end{center}

we conclude that $\log(I(M))=-2\log(\det(\triangle_{0})).$ Theorem \ref{IZ} is
proved. $\blacksquare.$

\subsection{Some Invariants of the Short Term Asymptotic Expansion of the Heat Kernel}

From the well know fact that for small $t$ we have $Tr(\exp(-t\triangle
_{0}))=\frac{vol(N)}{t^{n}}+\frac{k(g)}{t^{n-1}}+..+a_{0}+h(t,\tau
,\overline{\tau}),$ we will deduce that:

\begin{theorem}
\label{And1}Suppose that M is a CY manifold and g is a CY metric with a fixed
class of cohomology, then the coefficients \ a$_{k}$ for k=0,.,n \ in the
expression \ defined above are constant which depends only on the CY manifolds
and the fixed class of cohomology of the CY metric
\end{theorem}

\subsubsection{Idea of the Proof}

\textbf{\ }We know that the moduli space of CY metrics g with fix class of
cohomology is the same as the moduli space of complex structures. This follows
directly from the uniqueness and existence of the solution of the Calabi
problem. See \cite{Yau}. From here and results of Kodaira it follows that
$Tr(\exp(-t\triangle_{0}))$ is a smooth function with respect the coordinates
$\tau=(\tau_{1},.,\tau_{N})$ of the Kuranishi space $\mathcal{K}$(M). If we
prove that

\begin{center}
$\underset{t\rightarrow0}{\lim}\frac{\partial}{\partial\tau}\left(
Tr(\exp(-t\triangle_{0}))\right)  =\underset{t\rightarrow0}{\lim}%
\frac{\partial}{\partial\tau}\left(  \frac{vol(N)}{t^{n}}+\frac{k(g)}{t^{n-1}%
}+..+a_{0}+h(t,\tau,\overline{\tau})\right)  =c<\infty,$
\end{center}

then this implies that $a_{k}$ for $k=0,..,n=\dim_{\mathbb{C}}M$ are constants
on the moduli space.

\textbf{PROOF:} Let $F^{\prime}(1,\phi)$ be the operator defined in Definition
\ref{Hilb} and let $\sum a_{ii}$ be its trace, then we have:

\begin{lemma}
\label{And3} The following formula is true:
\end{lemma}

\begin{center}
$\frac{\partial}{\partial\tau}\left(  Tr(\exp\left(  -t\triangle_{\tau
,0}\right)  )\right)  |_{\tau=0}=t\sum_{i=1}^{\infty}\exp(-t\lambda
_{i})\lambda_{i}a_{ii}<\infty,$
\end{center}

\textit{for all }$t\geq0.$ $\lambda_{i}$ \textit{are eigen values of the
Laplacian} $\triangle_{0}$ on M$_{0}$ \textit{and }$.$

\textbf{PROOF:\ }According to Theorem \ref{t11} the following formula is true:

\begin{center}
$\frac{\partial}{\partial\tau}\left(  Tr(\exp\left(  -t\triangle_{\tau,0}%
^{"}\right)  \right)  |_{\tau=0}=-tTr\left(  \frac{d}{dt}\left(  \exp\left(
-t\triangle_{0}^{"}\right)  \right)  \circ\overline{\partial}^{-1}\circ
F^{\prime}(1,\phi)\circ\partial\right)  .$
\end{center}

According to Theorem \ref{holder} we have

\begin{center}
$Tr\left(  \frac{d}{dt}\left(  \exp\left(  -t\triangle_{0}^{"}\right)
\right)  \circ\overline{\partial}^{-1}\circ F^{\prime}(1,\phi)\circ
\partial\right)  =-t\sum_{i=1}^{\infty}\lambda_{i}\exp(-t\lambda_{i}%
)a_{ii}<\infty.$
\end{center}

Lemma \ref{And3} is proved. $\blacksquare.$

Since $\left\{  \lambda_{n}\right\}  $ for $n\geq1$ are the eigen values of
the Laplacian $\triangle_{0}$ then

\begin{center}
$\underset{n\rightarrow\infty}{\lim}\lambda_{n}\left(  (n^{\frac{2}{\dim M}%
})^{-1}\right)  =C>0,$
\end{center}

where dim$_{\mathbb{C}}$M is the complex dimension of M. From here and the
fact that the kernel of the operator $\Phi$ is a matrix with C$^{\infty}$
coefficients we derive that $\sum_{i=1}^{\infty}\lambda_{i}a_{ii}<\infty.$ So
Lemma\ \ref{And3} implies that

\begin{center}
$\underset{t\rightarrow0}{\lim}\frac{\partial}{\partial\tau}\left(
Tr(\exp\left(  -t\triangle_{0}\right)  \right)  =\underset{t\rightarrow0}%
{\lim}t\sum_{i=1}^{\infty}\lambda_{i}\exp(-t\lambda_{i})a_{ii}=0.$
\end{center}

This will imply that $\frac{\partial}{\partial\tau}\left(  a_{k}\right)  =0, $
for $k=0,.,n.$ Since $Tr(\exp(-t\triangle_{0}))$ is a real function when
$t\in\mathbb{R}$ and $t>0.$ So $\frac{\partial}{\partial\tau}\left(
a_{k}\right)  =0$ implies Theorem \ref{And1} directly. Theorem \ref{And1} is
proved. $\blacksquare.$

\section{The Analytic Torsion on CY Manifolds is Bounded.}

\begin{theorem}
\label{LY}Let M be any CY manifold, then 0$\leq\det(\triangle_{q})\leq C_{q}$
for $0\leq q\leq n=\dim M.$
\end{theorem}

\subsection{Outline of the Proof that of the Ray Singer Torsion on CY
Manifolds is bounded}

Theorem \ref{Bochner} implies that in order to prove Theorem \ref{LY} it is
enouph to bound $\det(\triangle_{0}).$ The bound of $\det(\triangle_{0})$ is
based on the following expression for the zeta function of the Laplacian
acting on functions:

\begin{center}
$\zeta_{0}(s)=\frac{1}{\Gamma(s)}\int_{0}^{\infty}\left(  Tr(\exp
(-t\triangle_{0}\right)  )t^{s-1}dt=b_{0}+b_{1}s+O(s^{2}).$
\end{center}

From the definition of $\det(\triangle_{0})$ it follows that $\det
(\triangle_{0})=\exp(-b_{1}).$ So if we bound $b_{1}$ Theorem \ref{LY} will be
proved. The bound $b_{1}$ is based on two facts. The first one is the
following asymptotic expansion of the Tr(exp($-t\triangle_{0}):$

\begin{center}
$Tr(\exp(-t\triangle_{0}))=\left(  \sum_{k=0}^{n}\frac{a_{k}}{t^{k}}\right)
+O(t)=\frac{vol(N)}{t^{n}}+\frac{k(g)}{t^{n-1}}+..+a_{0}+O(t),$
\end{center}

where n=dim$_{\mathbb{C}}M$ and k(g) is the scalar curvature of g.\textit{\ }%
See\textit{\ \cite{Roe}} on page 79. The second one is the explicit formula
for $b_{1}$ in \cite{ABKS}:

\begin{center}
$b_{1}=\gamma a_{0}+\sum_{k=1}^{n}\frac{a_{k}}{k}+\int_{0}^{1}\left(
Tr(\exp(-t\triangle_{0}))-\sum_{k=0}^{n}\frac{a_{k}}{t^{k}}\right)  \frac
{dt}{t}+\int_{1}^{\infty}Tr(\exp(-t\triangle_{0}))\frac{dt}{t},$
\end{center}

where $\gamma$ is the Euler constant. We will show that $b_{1}=C+\psi
(t,\tau,\overline{\tau}),$ where $\psi(t,\tau,\overline{\tau})\geq0.$ From
here we will obtain that $\det(\Delta_{0})\leq\exp(-C).$

\subsection{Proof of that the Analytic Torsion is Bounded}

\begin{remark}
\label{Y}From now on we will consider the following situation: \textit{We will
restrict our function }$h(t,\tau,\overline{\tau})$ \textit{on an one
dimensional disk in the Teichm\"{u}ller space $\mathcal{T}$}(M) \textit{of M
and this disk is defined as follows: }$\phi\in\mathbb{H}^{1}(M,T^{1,0})$
\textit{and let }
\end{remark}

\begin{center}
$\phi(\tau):=\phi\tau+\frac{1}{2}\overline{\partial}^{\ast}G[\phi(\tau
),\phi(\tau)].$
\end{center}

\textit{then we know that the Beltrami differential }$\phi\left(  \tau\right)
$\textit{\ \ is well defined }$C^{\infty}$\textit{\ section of }$C^{\infty
}(M_{0},Hom(\Omega^{1,0},\Omega^{0,1}))$ \textit{in a small disk for }%
$|\tau|<\varepsilon$ \textit{in the Teichm\"{u}ller space $\mathcal{T} $}(M).

\begin{theorem}
\label{And2}$-\log(\det(\Delta_{0}))=b_{1}(\tau,\overline{\tau})=C+\psi
(\tau,\overline{\tau}),$ where $C$ is a constant and $\psi(\tau,\overline
{\tau})\geq0.$
\end{theorem}

\textbf{PROOF:} Let us define $\psi_{1}(\tau,\overline{\tau}):=\frac
{\overline{\partial}^{2}}{\partial\tau\text{ }\overline{\partial\tau}}%
b_{1}(\tau,\overline{\tau}).$ Theorem \ref{var} implies that $\psi_{1}%
(\tau,\overline{\tau})\geq0.$ Let us define

\begin{center}
$\psi(\tau,\overline{\tau}):=\frac{1}{2\pi\sqrt{-1}}\iint_{|w-\tau|\leq1}%
\psi_{1}(\tau,\overline{\tau})G(\tau,w)d(w)\wedge\overline{d(w)},$
\end{center}

where $G(\tau,w)=-\log|w-\tau|$ is the Green kernel of the Laplacian
$\frac{\overline{\partial}^{2}}{\partial\tau\text{ }\overline{\partial\tau}}.$
Clearly since $\psi_{1}(\tau,\overline{\tau})\geq0$ and $G(\tau,w)\geq0$ for
$|\tau-w|\leq1$ we can conclude that $\psi(\tau,\overline{\tau})\geq0.$ From
the definition of the Green kernel we obtain that $\frac{\overline{\partial
}^{2}}{\partial\tau\text{ }\overline{\partial\tau}}\psi(\tau,\overline{\tau
})=\frac{\overline{\partial}^{2}}{\partial\tau\text{ }\overline{\partial\tau}%
}b_{1}(\tau,\overline{\tau})=\psi_{1}(\tau,\overline{\tau}).$ This fact
implies that b$_{1}(\tau,\overline{\tau})=$ $\psi(\tau,\overline{\tau}%
)+g(\tau)+\overline{g(\tau)},$ where $g(\tau)$ is a complex analytic function
in the disk $D$ defined in Remark \ref{Y}.

\begin{lemma}
\label{1000}$g(\tau)=const.$
\end{lemma}

\textbf{PROOF: }According to \cite{Roe} we have the following expression for
$Tr(\exp(-t\triangle_{0})):$

\begin{center}
$Tr(\exp(-t\triangle_{0}))=\left(  \sum_{k=0}^{n}\frac{a_{k}}{t^{k}}\right)
+O(t)=\frac{vol(N)}{t^{n}}+\frac{k(g)}{t^{n-1}}+..+a_{0}+h(t,\tau
,\overline{\tau}).$
\end{center}

According to \cite{ABKS} we have the following formula for $b_{1}%
(\tau,\overline{\tau}):$

\begin{center}
$\left(  \frac{d}{ds}\zeta_{0}(s)\right)  |_{s=0}=b_{1}(\tau,\overline{\tau
})=\gamma a_{0}+\sum_{k=1}^{n}\frac{a_{k}}{k}+\int_{0}^{1}\left(
Tr(\exp(-t\triangle_{0}))-\sum_{k=0}^{n}\frac{a_{k}}{t^{k}}\right)  \frac
{dt}{t}+\int_{1}^{\infty}Tr(\exp(-t\triangle_{0}))\frac{dt}{t}.$
\end{center}

Next we will compute $\frac{d}{d\tau}\left(  \frac{d}{ds}\zeta_{0}(s)\right)
. $ We have for large $s$ the following formula for

\begin{center}
$\zeta_{0}(s):=\frac{1}{\Gamma(s)}\int_{0}^{1}\left(  a_{0}+\sum_{k=1}%
^{n}\frac{a_{k}}{t^{k}}\right)  t^{s-1}dt+$

$\frac{1}{\Gamma(s)}\left(  \int_{0}^{1}\left(  Tr(\exp(-t\triangle_{0}%
))-\sum_{k=0}^{n}\frac{a_{k}}{t^{k}}\right)  t^{s}\frac{dt}{t}+\int
_{1}^{\infty}Tr(\exp(-t\triangle_{0}))t^{s}\frac{dt}{t}\right)  .$
\end{center}

So direct computations and Theorem \ref{And1} show that we have

\begin{center}
$\frac{d}{d\tau}\zeta_{0}(s)=\frac{d}{d\tau}\frac{1}{\Gamma(s)}\left(
\int_{0}^{1}Tr(\exp(-t\triangle_{0}))t^{s}\frac{dt}{t}+\int_{1}^{\infty
}Tr(\exp(-t\triangle_{0}))t^{s}\frac{dt}{t}\right)  =$

$\frac{d}{d\tau}\frac{1}{\Gamma(s)}\int_{0}^{\infty}Tr(\exp(-t\triangle
_{0}))t^{s}\frac{dt}{t}=\frac{1}{\Gamma(s)}\int_{0}^{\infty}Tr(\frac{d}{d\tau
}\exp(-t\triangle_{0}))t^{s}\frac{dt}{t}$
\end{center}

According to Lemma \ref{And3}:

\begin{center}
$\frac{d}{d\tau}\left(  Tr(\exp\left(  -t\triangle_{\tau}\right)  )\right)
=t\sum_{i=1}^{\infty}\lambda_{i}\exp(-t\lambda_{i})a_{ii}<\infty.$
\end{center}

where $Tr(F(1,\phi)=\sum_{i=1}^{\infty}a_{ii}.$ Combining these facts we
conclude that

\begin{center}
$\frac{d}{d\tau}b_{1}\left(  \tau,\overline{\tau})\right)  =\frac{d}%
{ds}\left(  \zeta_{\Delta_{0}}(s)\right)  |_{s=0}=\frac{d}{ds}\left(  \frac
{1}{\Gamma(s)}\int_{0}^{\infty}\left(  Tr(\frac{d}{d\tau}\exp(-t\triangle
_{0}))\right)  t^{s}\frac{dt}{t}\right)  |_{s=0}=$

$\frac{d}{ds}\left(  \frac{1}{\Gamma(s)}\int_{0}^{\infty}\left(  t\sum
_{i=1}^{\infty}\lambda_{i}\exp(-t\lambda_{i})a_{ii}\right)  t^{s}\frac{dt}%
{t}\right)  |_{s=0}=\frac{d}{ds}\left(  \sum_{i=1}^{\infty}\lambda_{i}%
^{s}a_{ii}\right)  |_{s=0}=\sum_{i=1}^{\infty}\log\left(  \lambda_{i}\right)
a_{ii}.$
\end{center}

Kodaira proved that the positive eigen values of the Laplacians $\overline
{\partial}_{\tau}\circ$ $\overline{\partial}_{\tau}$ depend on a $C^{\infty}$
manner in a small neighborhood of $\tau_{0}\in D.$ See \cite{KM}. From here
and the formula:

\begin{center}
$\frac{d}{d\tau}b_{1}\left(  \tau,\overline{\tau})\right)  =\frac{d}%
{ds}\left(  \zeta_{\Delta_{0}}(s)\right)  |_{s=0}=\frac{d}{d\tau}\left(
g(\tau)\right)  +\frac{d}{d\tau}\left(  \psi(t,\tau,\overline{\tau})\right)
=\sum_{i=1}^{\infty}\log(\lambda_{i})a_{ii}$
\end{center}

we can conclude that $\frac{d}{d\tau}\left(  g(\tau)\right)  =0.$ Lemma
\ref{1000} is proved. $\blacksquare.$ Lemma \ref{1000} implies Theorem
\ref{LY} directly. $\blacksquare.$

\end{document}